# Decomposable Branching Processes and Viral Marketing

Ranbir Dhounchak  
IEOR, IIT Bombay

Veeraruna Kavitha  
IEOR, IIT Bombay

**Abstract**

The decomposable branching processes are relatively less studied objects, particularly in the continuous time framework. In this paper, we consider various variants of decomposable continuous time branching processes. As usual practice in the theory of decomposable branching processes, we group various types into *irreducible classes*. These irreducible classes evolve according to the well-studied non-decomposable/irreducible branching processes. And we investigate the time evolution of the population of various classes when the process is initiated by the other class particle(s). We obtained class-wise extinction probability, and the time evolution of the population in the different classes. We then studied another peculiar type of decomposable branching process where any parent at the transition epoch either produces a random number of offspring, or its type gets changed (which may or may not be regarded as new offspring produced depending on the application). Such processes arise in modeling the content propagation of competing contents in online social networks. Here also, we obtain various performance measures. Additionally, we conjecture that the time evolution of the expected number of shares (different from the total progeny in irreducible branching processes) is given by the sum of two exponential curves corresponding to the two different classes.

## 1 Introduction

Earlier literature on studying the population dynamics considers deterministic models, such as differential equations, difference equations, etc. While these models being elegant in many aspects, they have serious drawbacks which make them inadequate to capture the critical components of the dynamics. To mention a few, the reproductive behavior of individuals is not deterministic as considered in deterministic models; it is rather stochastic. Further, there are other sources of randomness which are not incorporated in deterministic models such as the distribution of life span of individuals, the environment, etc. These shortcomings led to the introduction of stochastic models which incorporate critical features/characteristics of the population dynamics. And branching processes are considered to be a well-suited candidate for modeling the population dynamics. Branching processes are popular stochastic processes that describe the evolution of the population over generations. These are proven to be an adequate tool for modeling the dynamics of 'population' for various applications like tumor cells proliferation [8, 9], advertising over social network [10], etc. The key questions in the theory of branching processes include:

- What is the extinction probability?
- What is the growth rate of the population?
- What is the total progeny? etc.

Branching processes can be categorized on a number of factors, for example, the discrete and continuous time branching processes (classification by time), single type and multi-type branching processes, critical, super-critical or sub-critical, etc. Further, each category has subcategories giving rise to numerous variants of the branching processes.

In a multitype *continuous time branching process* (CTBP), a particle lives for an exponentially distributed random time. It produces a random number of offspring of various types independent of the other particles and then dies. And this continues. The underlying generator matrix, say $A$, plays a vital role in carrying out the analysis of CTBP. When the $e^{At}$ matrix is called positive regular[1], the underlying CTBP is classified as irreducible/non-decomposable. In this case, the branching processes exhibit a certain dichotomy (e.g., [2, 4]): a) either all the types survive together and grow exponentially (with time) with the same rate; b) or all the types get extinct after some time. Mainly, the largest eigenvalue of $A$, say $\alpha$, determines the growth rate, extinction probability, etc. The CTBP is called subcritical, critical and supercritical based on whether $\alpha < 1, \alpha = 1$ and $\alpha > 1$ respectively. When $\alpha \leq 1$, the population gets extinct with probability one. Whereas when $\alpha > 1$, the CTBP can survive on some sample paths. Further, on these sample paths, all types grow exponentially fast with the common rate $\alpha$ provided the CTBP is non-decomposable [4].

---

[1]Matrix $e^{At}$ is positive regular if each entry of it is strictly positive for some $t_0 > 0$.

When the process is such that the particles of certain types do not produce offspring of certain other types, we have a very different variety of branching process called as decomposable branching process. In this process, the types get partitioned into different classes, where the types across different classes may have different characteristics. These processes behave significantly different from non-decomposable processes. First and foremost, the dichotomy no longer holds, i.e., a particular class (a group of types) may thrive/survive whereas another gets extinct. Secondly, the types in different classes may have different growth rates, etc.

In this paper, we analyze various variants of decomposable branching processes in the continuous time framework. We first study the process where any irreducible class contains only one type of particles, which we refer to as scalar decomposable branching processes (SDCBP). We then consider the vector version, i.e., the case with two irreducible classes with each class having particles of multiple types. We analyze the growth rates of these irreducible classes on their respective survival paths. We provide the almost sure analysis of these processes on the viral paths using continuous time martingale theory. Next, we study a peculiar version of the branching processes where at each transition epoch the parent either produces offspring or gets changed to a different type. This is inspired by content propagation in *online social network* (see [3]). We call these as 'type-changing decomposable branching processes.' Here we study two types of total progeny: a) one which counts the type changes as new offspring; b) one which does not count type changes as new offspring. We conjecture the expected time evolution of the total progeny by solving fixed point equations of certain integral operators in infinite dimensional Banach spaces.

## 1.1 Related literature

There is a vast literature that studies non-decomposable branching processes (e.g., [2, 4, 7], etc). For example, Wang et al. [7] computed the probability generating function of the total progeny of multi-type irreducible branching random walk. They showed that asymptotically the distribution of the total progeny satisfies a fixed point equation involving probability generating function of offspring. One can find a lot more work in this direction, and the references for the same can be found in [2, 4, 7] etc.

The decomposable processes are relatively less studied in comparison with the non-decomposable objects. There are strands of literature (e.g., [6, 5]) that study the discrete time decomposable branching processes. But to the best of our knowledge, there is no dedicated literature on decomposable continuous time branching process (DCTBPs). Some analysis could probably be derived using the discrete-time results. However, there are many more questions that need to be answered directly for continuous time versions, and we consider the same. We also consider more questions related to 'type-changing' decomposable branching processes.

Kesten et al. [6] studied the multitype decomposable branching process in the discrete time framework. Their study mainly focused on investigating the growth rates of the particles belonging to different classes. Hautphenne [5] investigated the class-wise extinction probabilities, where the extinction of a class is shown to be the minimal non-negative solution of the extinction probability equation but with added constraints. Again, no work is done in investigating the growth rates of various types in the *continuous time framework*. The continuous time version is useful in studying many practical problems such as cancer biology, viral marketing problem[3]. In these problems, it becomes important to know the following: the growth rates of different classes, when do the particles of one specific class explode while the others get extinct? etc. And all of these measures require the theory of continuous time decomposable branching processes.

The organization of this paper is as follows. In section 2, we study the scalar decomposable branching processes, while section 3 considers the vector version. Type-changing processes are considered in section 4. The appendix contains the proofs.

## 2 Scalar decomposable branching processes

We consider a decomposable continuous time branching process with $n$ types of particles. A particle of type-$i$ lives for an exponentially distributed time with parameter $\lambda$ and produces the offspring of types $j \geq i$ only. Each parent independently (of the others) produces a random number of offspring and according to an identical distribution when the parents are of the same type. It is a common practice in the theory of decomposable branching processes to group various types into *irreducible classes*. Any irreducible class, say $\mathscr{C}$, consists of all those types such that the class forms a non-decomposable/irreducible branching process when the process is initiated by particles within $\mathscr{C}$.

In this section, we consider the scenario in which each class has exactly one type. We refer to this process as *scalar decomposable continuous time branching processes* (SDCBP). In the later sections, we generalize it to classes having two or more types of particles.

The mathematical framework of SDCBP is as follows. Let $\mathbf{X}(t) = \{X_1(t), X_2(t), \cdots, X_n(t)\}$ denote the population vector where $X_i(t)$ represents the number of type-$i$ particles at time $t$. Let $\zeta_{ij}$ be the number of type-$j$ offspring produced by a type-$i$ particle (note that $\zeta_{ij} = 0$ for $j < i$). The probability generating function of offspring of a particle of type-$i$, say $h_i(\mathbf{s})$, is given as:

$$h_i(\mathbf{s}) := \sum_{k_i, \cdots, k_n} P(\zeta_{ii} = k_i, \zeta_{ii+1} = k_{i+1}, \cdots, \zeta_{in} = k_n) \Pi_{j \geq i} s_j^{k_j}, \quad \mathbf{s} = (s_1, \cdots, s_n) \in [0,1]^n. \tag{1}$$

The *generator matrix* plays an instrumental role in carrying out the study of the branching process. We now obtain the generator matrix, say $B$, of the underlying SDCBP. Denote by $B_{ij}$ the $(i,j)$-th entry of $B$ which is $\alpha_{ij} - \delta_{ij}$ where $\alpha_{ij} := \lambda \frac{\delta h_i(\mathbf{s})}{\delta s_j}\big|_{\mathbf{s}=\mathbf{1}}$ and $\delta_{ij} = 1$ when $i = j$ and 0 else (as in [4]). In our case, the generator matrix $B$ is:

$$B := B_n = \begin{bmatrix} \alpha_{11} - \lambda & \alpha_{12} & \cdots & \alpha_{1n-1} & \alpha_{1n} \\ 0 & \alpha_{22} - \lambda & \cdots & \alpha_{2n-1} & \alpha_{2n} \\ \vdots & \vdots & \ddots & \vdots & \vdots \\ 0 & 0 & \cdots & 0 & \alpha_{nn} - \lambda \end{bmatrix} = \begin{bmatrix} \alpha_1 & \alpha_{12} & \cdots & \alpha_{1n-1} & \alpha_{1n} \\ 0 & \alpha_2 & \cdots & \alpha_{2n-1} & \alpha_{2n} \\ \vdots & \vdots & \ddots & \vdots & \vdots \\ 0 & 0 & \cdots & 0 & \alpha_n \end{bmatrix}; \tag{2}$$

with $\alpha_i := \alpha_{ii} - \lambda \ \forall i$. We immediately see that the matrix $B$ is *not* positive regular[2], so the underlying process is decomposable.

There are many ways to characterize decomposable branching processes, and we characterize it by the underlying generator matrix. A continuous type branching process is called non-decomposable/irreducible if the matrix $exp(B)$[3] has *all entries* strictly positive and decomposable otherwise (e.g., see [5]). With this understanding, we are now ready to address various branching process related questions in the context of SDCBP.

## 2.1 Extinction probability

We begin with investigating the extinction probability of SDCBP. The notion of extinction probability is not the same as in irreducible branching processes. Here we study type-wise extinction, unlike the irreducible branching processes. Observe that a type-1 particle can be produced only by itself, and consequently, its extinction does not depend on the population of other types. In general, the extinction probability of type-$i$ particles depends only on its predecessor types including itself, i.e., type-1,type-2,$\cdots$, type-$i$; which can also be seen through the generator matrix. In other words, type-$i$ particles get extinct when all the particles of type-$k$ with $k \leq i$ get extinct. We adopt the convention that the extinction probability of type-$i$ particles is always *one* when the process is initiated by any type $j \geq i + 1$ particle; as none of the particles of type $j \geq i + 1$ produce offspring of type-$i$. Define the following

$$q_{ki} := P\big[X_j(t) = 0 \text{ for all } j \leq i \text{ for some } t \mid \mathbf{X}(0) = \mathbf{e}_k\big], \tag{3}$$

where $\mathbf{e}_k$ is a $n$-dimensional unit vector with 1 only at $k$-th place. Observe that the above definition of $q_{ki}$ is in accordance with the said notion of the extinction of the type-$i$ population. Note that $q_{ki} = 1$ when $k > i$. And when $k \leq i$, we obtain these extinction probabilities by solving the fixed point equations (see, e.g. [3]) via conditioning on the events of the first transition

$$q_{ki} = h_k\left(\sum_{j=1}^{k-1} \mathbf{e}_j + \sum_{j=k}^{i} q_{ji} \mathbf{e}_j + \sum_{j=i+1}^{n} \mathbf{e}_j\right). \tag{4}$$

This is due the following. At the first transition, the type-$k$ particle produces the offspring of type-$j$ (recall $j \geq k$) based on which we have the following two scenario:

- when $j > i$, the extinction probability of type-$i$ is one;

---
[2] Matrix $B$ is called positive regular if there exists a $n$ such that the matrix $B^n$ has all strict positive entries, i.e., $B_{ij}^n > 0 \ \forall i, j$.

[3] $exp(B) = e^B = I + B + \frac{B^2}{2!} + \frac{B^3}{3!} + \cdots$

- when $j \le i$, $\zeta_{kj}$ offsprings are produced, the extinction probability becomes $q_{ji}^{\zeta_{kj}} \times P(\zeta_{kj})$. By conditioning on number of offspring produced, we have the expression for $q_{ki}$ as given in equation (4).

We illustrate how a few types can get extinct, and a few others can thrive as follows. Consider a decomposable branching process initiated with one type-1 particle. Say it produces the offspring of all but its own type upon its death. Thus, in this event, the type-1 particle gets extinct. Whereas the particles of remaining types keep producing offspring in subsequent generations and their populations may explode.

## 2.2 Time evolution of the population

The type-$i$ population of the SDCBP, when initiated with the type-$i$ particles, evolves on its own while influencing the evolution of the successor/higher types ($> i$) particles. Further, the evolution of type-$i$ population is well-understood in the literature (e.g., [2, 4]). And we aim to study the evolution of the other types starting from a different type of particle.

Denote by $\mathscr{F}_{i,t} := \sigma\{X_i(u), u \le t\}$ the sigma algebra generated by type-$i$ population. The stochastic processes $\{X_i(t)e^{-\alpha_i t}; \mathscr{F}_{i,t}; t \ge 0\}$ is a *convergent martingale* when it starts from type-$i$ particles (see [4]). Further
$$\lim_{t\to\infty} X_i(t)e^{-\alpha_i t} = W_i; \quad \text{where } W_i \ge 0. \text{ Equivalently } X_i(t)e^{-\alpha_i t} \xrightarrow{a.s.} W_i \text{ as } t \to \infty.$$
The non-negative random variable $W_i$ satisfies the following: when $E[\zeta_{ii}\log\zeta_{ii}] < \infty\ \forall i$, the extinction probability of type-$i$ particles $q_{ii}$ (when started with self type) is given by $P(W_i = 0|\mathbf{X}(0) = \mathbf{e}_i) = q_{ii}$ and $E[W_i|\mathbf{X}(0) = \mathbf{e}_i] = 1$ (see [2, 4]).

To the best of our knowledge, the evolution of the population of type-$i+1,\cdots$, type-$n$ when initiated with type-$i$ particles is *not* studied in the continuous time framework. And we precisely investigate this evolution. It is easy to see that when the $\alpha := \max_i \alpha_i$ is *strictly less than one*, the whole process is sub-critical and $\mathbf{X}(t) = \mathbf{0}$ for some $t > 0$. On the other hand, when $\alpha > 0$, we have the super-critical process. We aim to derive the analysis of the super-critical decomposable branching process. At a more general sense, we study the evolution of type-$m$ population when started with one type-1 particle. The Theorem below provides the same:

**Theorem 1** *Let $(\Omega, \mathscr{F}_t, \mathbb{P})$ be the probability space and $\mathscr{F}_t := \sigma\{\mathbf{X}(u), u \le t\}$ be the natural filtration for the underlying branching process. Let $E[\zeta_{ij}] < \infty$ for all $i, j$. When started with a particle of type-1, we have:*

(i) (a) *The following stochastic process for any $2 \le m \le n$ is a martingale when $\alpha_1 \ne \alpha_2 \ne \cdots \ne \alpha_n$*

$$X_m(t)e^{-\alpha_m t} + \sum_{k=1}^{m-1}\sum_{j_k=j_{k-1}+1}^{m-1}\cdots\sum_{j_3=j_2+1}^{m-k+2}\sum_{j_2=j_1+1}^{m-k+1}\sum_{j_1=1}^{m-k}\frac{\alpha_{j_1 j_2}\alpha_{j_2 j_3}\cdots\alpha_{j_k m}X_{j_1}(t)e^{-\alpha_m t}}{(\alpha_m-\alpha_{j_1})(\alpha_m-\alpha_{j_2})\cdots(\alpha_m-\alpha_{j_k})}. \quad (5)$$

*An alternative representation of the above martingale is:* $\sum_{i=1}^{m} a_i^m X_i(t)e^{-\alpha_m t}$ *where $a_i^m = 1\forall\ i$, and for $i = 1,2\cdots,m-1$*

$$a_i^m = \sum_{k=1}^{m-i-1}\sum_{j_k=j_{k-1}+1}^{m-1}\cdots\sum_{j_1=i+1}^{m-k}\frac{\alpha_{ij_1}\alpha_{j_1 j_2}\cdots\alpha_{j_k m}}{(\alpha_m-\alpha_i)(\alpha_m-\alpha_{j_1})(\alpha_m-\alpha_{j_2})\cdots(\alpha_m-\alpha_{j_k})} + \frac{\alpha_{im}}{\alpha_m-\alpha_i}.$$

(b) $E[X_m(t)] = a_1^m e^{\alpha_m t} + \sum_{j=1}^{m-1} P_j^m e^{\alpha_j t}\ \forall t$, *where*

$$P_1^m = -a_1^m - \sum_{k=1}^{m-2}(-1)^k\sum_{j_1=2}^{m-k}\sum_{j_2=j_1+1}^{m-k+1}\cdots\sum_{j_k=j_{k-1}+1}^{m-1} a_1^{j_1} a_{j_1}^{j_2}\cdots a_{j_k}^m, \quad P_{m-1}^m = -a_1^{m-1} a_{m-1}^m; \text{ and}$$

$$P_j^m = \sum_{k=1}^{m-j-1}(-1)^k\sum_{j_1=j+1}^{m-k}\cdots\sum_{j_k=j_{k-1}+1}^{m-1} a_1^j a_j^{j_1} a_{j_1}^{j_2}\cdots a_{j_k}^m - a_1^j a_j^m; \text{ for } j = 2,3\cdots,m-2. \quad (6)$$

(ii) *Further, when $\alpha_1 < \alpha_2 < \cdots < \alpha_n$, we have a non negative martingale and hence it converges almost surely to a non-negative random variable, say $W_m$, which is integrable. Thus, as $t \to \infty$*

$$X_m(t)e^{-\alpha_m t} + \sum_{k=1}^{m-1}\sum_{j_1=1}^{m-k}\sum_{j_2=j_1+1}^{m-k+1}\sum_{j_3=j_2+1}^{m-k+2}\cdots\sum_{j_k=j_{k-1}+1}^{m-1}\frac{\alpha_{j_1 j_2}\alpha_{j_2 j_3}\cdots\alpha_{j_k m}X_{j_1}(t)e^{-\alpha_m t}}{(\alpha_m-\alpha_{j_1})(\alpha_m-\alpha_{j_2})\cdots(\alpha_m-\alpha_{j_k})} \xrightarrow{a.s.} W_m. \quad (7)$$

*In addition, we also have*

$$\left| X_m(t)e^{-\alpha_m t} - \sum_{i=1}^{m-1}\frac{P_i^m}{a_1^i} W_i e^{(\alpha_i-\alpha_m)t} - W_m \right| \xrightarrow{a.s.} 0; \text{ where } W_1, W_2, \cdots, W_m \text{ are non-negative random variables.}$$

**Proof 1** *The proof is given in Appendix.* ∎

We immediately see that after sufficiently large time $t$, $X_m(t)$ grow exponentially with rate $\alpha_m$ when the process is initiated with type-1 particle. In general, one can write the time evolution of type-$m$ particle when the process starts with type-$p$ particle with $m \geq p$ in a similar fashion.

We next show that the characteristics of the random variables $\{W_m\}_m$ in SDCBP are similar to that in the irreducible branching processes (e.g., see [4]).

**Theorem 2** *With $E[\zeta_{ij} \log \zeta_{ij}] < \infty \; \forall i, j$ and $\alpha_m > \alpha_{m-1}, \cdots, \alpha_1$, we have $P(W_j = 0 | \mathbf{X}_0 = \mathbf{e}_k) = q_{kj}$.*

**Proof 2** *The proof is given in the Appendix.* ∎

The above result is similar to that in the non-decomposable branching process. However, here we have class wise extinction probability. In other words, the random variable $W_j$ characterizes the extinction probability of type-$j$ particles.

## 3 Vector decomposable continuous time branching processes

In the previous section, we studied the decomposable branching processes where each irreducible class contains only one type. We now extend this model to include the cases where a class contains more than one type of population. We consider a decomposable *branching process* (BP) with $n + m$ types. And as before the population types can be partitioned into several classes such that all the types within a class form an irreducible BP.

Let $\{\mathbf{Z}(t) = \mathbf{X}(t), \mathbf{Y}(t), t \in [0, \infty)\}$ be an $(n + m)$−type decomposable BP where $\mathbf{X}(t) = \{X_1(t), \cdots, X_n(t)\}$ and $\mathbf{Y}(t) = \{Y_1(t), \cdots, Y_m(t)\}$ both form irreducible multitype BP on their own. Let us say that these types are partitioned into two classes, i.e., the types of $\mathbf{X}(t)$ belong to class $\mathbb{C}_1$ and that of $\mathbf{Y}(t)$ belong to class $\mathbb{C}_2$. A particle belonging to class $\mathbb{C}_1$ produces the offspring of the types in $\mathbb{C}_2$ along with the offspring of the types of its own class. While the particles of $\mathbb{C}_2$ class produce offspring of their own class only. In this way, the $\mathbb{C}_1$ population influences the $\mathbb{C}_2$ population, but not the other way around. Say $\tilde{\zeta}_{ij}$ be the number of independent and identically distributed (IID) type-$j$ offspring produced by a type-$i$ particle. Note that when $i \in \mathbb{C}_2, j \in \mathbb{C}_1$, then $\tilde{\zeta}_{ij} = 0$. We refer to this process as *vector decomposable continuous time branching process* (VDCBP).

The evolution of the process $\mathbf{X}(t)$ alone is well understood in literature (see [4]). Observe that $\mathbf{Y}(t)$ also forms an irreducible branching process when $\mathbf{X}(t)$ is absent. In other words, the evolution of the stochastic process $\mathbf{Y}(t)$ is well understood when the process is initiated by a particle of $\mathbb{C}_2$. However, the evolution of $\mathbf{Y}(t)$ in presence of $\mathbf{X}(t)$, i.e., when started with particle(s) of $\mathbb{C}_1$ is not known in the continuous time framework. Authors in [6] study the evolution of a similar decomposable branching process in the discrete time framework. However, *no attempt has been made to study the underlying martingales and associated results even for the discrete case*.

### 3.1 Probability generating function and the time evolution of the population

Let $\tilde{h}_i(\mathbf{s})$ represent the probability generation function (PGF) of number of offspring produced by a type-$i$ particle, where $\mathbf{s} = (s_1, s_2, \cdots, s_{m+n}) \in [0,1]^{n+m}$. With $\mathbf{K} := k_1, \cdots, k_n, \cdots k_{m+n}$ the PGF is

$$\tilde{h}_i(\mathbf{s}) = \sum_{\mathbf{K}} P(\tilde{\zeta}_{i1} = k_1, \cdots, \tilde{\zeta}_{in} = k_n, \cdots, \tilde{\zeta}_{in+m} = k_{n+m}) \Pi_{j=1}^{n+m} s_j^{k_j}.$$

Let $\mathbb{A}$ be the (embedded) generator matrix of the branching process $\{\mathbf{Z}(t)\}$. The $(i, j)$−th entry of $\mathbb{A}$ represents the expected number of offspring of type-$j$ produced by a particle of type-$i$:

$$\mathbb{A}_{i,j} = \frac{\partial \tilde{h}_i(\mathbf{s})}{\partial s_j}\bigg|_{\mathbf{s}=\mathbf{1}}; \quad \mathbf{1} = [1, 1, \cdots, 1] \text{ n+m-dimensional vector}.$$

In what follows, the generator matrix $\mathbb{A}$ has the following structure

$$\mathbb{A} = \begin{bmatrix} A_{11} & A_{12} \\ \tilde{0} & A_{22} \end{bmatrix}; \quad A_{11} \in \mathbb{R}^{n \times n}, \; A_{22} \in \mathbb{R}^{m \times m}, \quad A_{12} \in \mathbb{R}^{n \times m}, \; \tilde{0} \in \mathbb{R}^{m \times n} \tag{8}$$

where $A_{11}, A_{22}$ correspond to the transitions within $\mathbb{C}_1$ and $\mathbb{C}_2$ respectively, while $A_{12}$ corresponds to the transitions from $\mathbb{C}_1$ to $\mathbb{C}_2$. The matrices $A_{11}, A_{22}$ are positive regular as the processes $\mathbf{X}(t)$ and $\mathbf{Y}(t)$ are irreducible when started with their own type particles. The largest eigenvalue of $e^{A_{ii}}$ determines the growth rates of $\mathbb{C}_i$ class population ($i = 1, 2$). And using the Perron Frobenius theory of positive regular matrices, there exists a positive real eigenvalue, say $e^{\alpha_i}$ the Perron root, of the matrix $e^{A_{ii}}$ such that: a) its algebraic and geometric multiplicities both equal to one; and b) it dominates other eigenvalues absolutely. Further, the left and right eigenvectors corresponding to the said Perron root can be taken in a way that each component of the both the vectors is positive and the inner product of both vectors is one. In other words, with $\xi_i^L$ and $\xi_i^R$ as the left eigenvector and right eigenvector of the matrix $e^{A_{ii}}$ corresponding to the Perron root $e^{\alpha_i}$, we have[4]

$$\xi_i^R . \xi_i^L = 1 \text{ and } \xi_i^L > \mathbf{0}, \; \xi_i^R > \mathbf{0}.$$

It is well-known that the process $\{\xi_1^R . \mathbf{X}(t) e^{-\alpha_1 t}\}$ is a martingale under the natural filtration, and $E[\mathbf{X}(t)] = e^{\alpha_1 t} \xi_1^L$ [4]. Also, the process $\{\xi_2^R . \mathbf{Y}(t) e^{-\alpha_2 t}\}$ is a martingale provided that the progenitor belongs to class $\mathbb{C}_2$ under the natural filtration. Further, both of these martingales converge to non-negative random variables. As before, we investigate the evolution of $\mathbb{C}_2$ particles when the process starts with a particle of the class $\mathbb{C}_1$. We have the following Theorem describing the same:

**Theorem 3** *Let $(\Omega, \mathbb{F}_t, \mathbb{P})$ be the probability space and $\{\mathbb{F}_t = \sigma(\mathbf{Z}(\omega, t), t \geq 0)\}$ be the natural filtration for the process $\{\mathbf{Z}(t)\}_t$. When the process starts with a type-$j$ particle of $\mathbb{C}_1$, we have the following result.*

  i. *The stochastic process $\{e^{-\alpha_2 t} \mathbf{Y}(t) \xi_2^R + e^{-\alpha_2 t} \mathbf{X}(t)(\alpha_2 I - A_{11})^{-1} A_{12} \xi_2^R; \; \mathbb{F}_t; \; t \geq 0\}$ becomes a martingale provided $\alpha_1 \neq \alpha_2$. Further, for any $t > 0$*

$$E\left[e^{-\alpha_2 t} \mathbf{Y}(t) \xi_2^R + e^{-\alpha_2 t} \mathbf{X}(t)(\alpha_2 I - A_{11})^{-1} A_{12} \xi_2^R\right] = \mathbf{e}_j (\alpha_2 I - A_{11})^{-1} A_{12} \xi_2^R. \tag{9}$$

  ii. *When $\alpha_1 < \alpha_2$, the martingale converges and hence we have the following as $t \to \infty$*

$$e^{-\alpha_2 t} \mathbf{X}(t) \xrightarrow{a.s.} \mathbf{0}_{n \times 1} \text{ and } \left| e^{-\alpha_2 t} \mathbf{Y}(t) \xi_2^R + e^{-\alpha_2 t} \mathbf{X}(t)(\alpha_2 I - A_{11})^{-1} A_{12} \xi_2^R \right| \xrightarrow{a.s.} \tilde{W}_{12}$$

  *where $\tilde{W}_{12}$ is a integrable random variable and $\mathbf{0}_{n \times 1}$ n-dimensional null column vector. As before, we also have $\left| e^{-\alpha_2 t} \mathbf{Y}(t) \xi_2^R + e^{-\alpha_2 t} W_1 e^{\alpha_1 t} \xi_1^R (\alpha_2 I - A_{11})^{-1} A_{12} \xi_2^R - \tilde{W}_{12} \right| \xrightarrow{a.s.} 0.$*

***Proof 3*** *The proof is given in Appendix.* ∎

Using equation 9, we have

$$\begin{aligned}
E\left[e^{-\alpha_2 t} \mathbf{Y}(t) \xi_2^R\right] + E\left[e^{-\alpha_2 t} \mathbf{X}(t)(\alpha_2 I - A_{11})^{-1} A_{12} \xi_2^R\right] &= \mathbf{e}_j (\alpha_2 I - A_{11})^{-1} A_{12} \xi_2^R; \text{ on multiplying with } \xi_2^L \\
e^{-\alpha_2 t} E\left[\mathbf{Y}(t)\right] \xi_2^R . \xi_2^L + E\left[e^{-\alpha_2 t} \mathbf{X}(t)(\alpha_2 I - A_{11})^{-1} A_{12}\right] \xi_2^R . \xi_2^L &= \mathbf{e}_j (\alpha_2 I - A_{11})^{-1} A_{12} \xi_2^R . \xi_2^L \\
E\left[\mathbf{Y}(t)\right] + E\left[\mathbf{X}(t)(\alpha_2 I - A_{11})^{-1} A_{12}\right] &= e^{\alpha_2 t} \mathbf{e}_j (\alpha_2 I - A_{11})^{-1} A_{12}; \quad \because \xi_2^R . \xi_2^L = 1 \\
E\left[\mathbf{Y}(t)\right] &= e^{\alpha_2 t} \mathbf{e}_j (\alpha_2 I - A_{11})^{-1} A_{12} - E\left[\mathbf{X}(t)\right](\alpha_2 I - A_{11})^{-1} A_{12}
\end{aligned}$$

As $E\left[\mathbf{X}(t)\right] = e^{\alpha_1 t} \xi_1^L$, one can write $E\left[\mathbf{Y}(t)\right] = e^{\alpha_2 t} \mathbf{e}_j (\alpha_2 I - A_{11})^{-1} A_{12} - e^{\alpha_1 t} \xi_1^L (\alpha_2 I - A_{11})^{-1} A_{12}$. With $E\left[\mathbf{X}(l, t)\right]$ as the $l$–th component of $E\left[\mathbf{X}(t)\right]$, then we have

$$E\left[\mathbf{Y}(l, t)\right] = \hat{h}_l e^{\alpha_2 t} - \hat{o}_l e^{\alpha_1 t} \tag{10}$$

$\hat{h}_l$ and $\hat{o}_l$ are the $l$-th component of the vectors $\mathbf{e}_j (\alpha_2 I - A_{11})^{-1} A_{12}$ and $\xi_1^L (\alpha_2 I - A_{11})^{-1} A_{12}$ respectively.

Thus, we see that expected value of the class $\mathbb{C}_2$ population when the process starts with the class $\mathbb{C}_1$ type particle is given by the sum of two exponential curves corresponding the classes $\mathbb{C}_1$ and $\mathbb{C}_2$.

**Extinction probability and limit:** The limiting behaviour of particles belonging $\mathbb{C}_i$, $i = 1, 2$ when the process starts with the particle from class (i.e., progenitor belongs to $\mathbb{C}_i$) is well-known (e.g., see [4]). Basically,

$$\zeta_1^R . \mathbf{X}(t) e^{-\alpha_1 t} \xrightarrow{a.s.} \tilde{W}_1; \quad \zeta_2^R . \mathbf{Y}(t) e^{-\alpha_2 t} \xrightarrow{a.s.} \tilde{W}_2 \tag{11}$$

---
[4]Vector inequality $\mathbf{V} > \mathbf{U}$ means that $v_j > u_j \; \forall \; j$; where $v_j, u_j$ are the $j$–th components of the vectors $\mathbf{V}$ and $\mathbf{U}$.

where $\tilde{W}_i$ is a non-negative random variable. Further, the random variable $\tilde{W}_i$ characterizes the extinction probability of $\mathbb{C}_i$ particles. Denote by $q_1^j$ the extinction probability of $\mathbb{C}_1$ particles when the progenitor is of type-$j$ from $\mathbb{C}_1$, i.e., $q_1^j := P\left(\mathbf{X}(t) = 0 \text{ for some } t > 0 \middle| \mathbf{X}(0) = \mathbf{e}_j, \mathbf{Y}(0) = \mathbf{0}\right)$. When $E[\tilde{\zeta}_{i,j} \log \tilde{\zeta}_{i,j}] < \infty \; \forall \; i, j$, then it is known that [4]

$$q_1^j = P\left(\tilde{W}_1 = 0 \middle| \mathbf{X}(0) = \mathbf{e}_j, \mathbf{Y}(0) = \mathbf{0}\right). \tag{12}$$

Similarly, the extinction probability of the class $\mathbb{C}_2$, $q_2^j := P\left(\mathbf{Y}(t) = 0 \text{ for some } t > 0 \middle| \mathbf{X}(0) = \mathbf{0}, \mathbf{Y}(0) = \mathbf{e}_j\right)$, satisfies $q_2^j := P\left(\tilde{W}_2 = 0 \middle| \mathbf{X}(0) = \mathbf{0}, \mathbf{Y}(0) = \mathbf{e}_j\right)$.

We are now left with the computation of the extinction probability of $\mathbb{C}_2$ particles when the process is initiated with a class $\mathbb{C}_1$ particle. The extinction in this scenario implies the full extinction, i.e., particles of both the classes get extinct. Define $q^{1j}$ as below:

$$q^{1j} := P\left(\mathbf{Y}(t) = 0 \text{ for some } t > 0 \middle| \mathbf{X}(0) = \mathbf{e}_j, \mathbf{Y}(0) = \mathbf{0}\right) = P\left(\mathbf{X}(t) = 0, \mathbf{Y}(t) = 0 \text{ for some } t > 0 \middle| \mathbf{X}(0) = \mathbf{e}_j, \mathbf{Y}(0) = \mathbf{0}\right)$$

Here also it is easy to see that $q^{1j} = P\left(\tilde{W}_{12} = 0 \middle| \mathbf{X}(0) = \mathbf{e}_j\right)$ provided that $E[\tilde{\zeta}_{i,j} \log \tilde{\zeta}_{i,j}] < \infty \; \forall \; i, j$ (similar to Theorem 2).

## 4 Type-changing vector decomposable branching processes

We now study a slightly different version of the above vector decomposable branching process. The difference lies in the transition/reproduction epoch: any parent at the transition epoch either produces a random number of offspring (as before) *or* its type gets changed (one of the two events takes place). And then it dies. We consider a decomposable branching process consists of two irreducible classes, namely mixed ($\mathcal{M}_x$) class and exclusive ($\mathcal{E}_x$) class. Particles of $\mathcal{M}_x$ class produce particles of $\mathcal{M}_x$ class as well as that of $\mathcal{E}_x$ class. While particles of $\mathcal{E}_x$ produces particles of $\mathcal{E}_x$ class only. In one of the two events at a transition epoch, a particle of $l \in \mathcal{M}_x$ type wakes up after exponentially distributed time with parameter $v$, i.e., $exp(v)$ and produces a random number of offspring. Whereas in the other event, the type of the particle gets changed after $exp(\lambda)$ time. It is easy to see that probability of the former event is $1 - \theta$ and that of the latter is $\theta$ where $\theta = \lambda/(\lambda + v)$. We refer to this process briefly as *type-changing vector decomposable branching process* (TC-VDBP).

Let $m_{l,k}$ represent the expected number of offspring of type $k$ produced by a parent of type $l$, where $l$ and $k$ can be of $\mathcal{E}_x$ class or $\mathcal{M}_x$ class. And $a_{l,k}$ is the probability that a $l \in \mathcal{M}_x$ particle gets converted to a $k \in \mathcal{M}_x$ particle. In a similar way, type change transitions are allowed within $\mathcal{E}_x$ class, however, there are no type changes possible from one class to another. The generator matrix of such a TC-VDBP has the following structure:

$$\begin{bmatrix} A_{mx} & A_{mx,ex} \\ \mathbf{0} & A_{ex} \end{bmatrix}, \tag{13}$$

where $A_{mx}$ represents all the transitions between types belonging to class $\mathcal{M}_x$, $A_{ex}$ represents all the transitions between types belonging to class $\mathcal{E}_x$, while $A_{mx,ex}$ represents the transition between $\mathcal{M}_x$ and $\mathcal{E}_x$ (offspring of class $\mathcal{E}_x$ produced by class $\mathcal{M}_x$). With the description as above for example:

$$A_{mx} := \begin{bmatrix} \theta a_{1,1} + (1-\theta) m_{1,1} & \theta a_{1,2} + (1-\theta) m_{1,2} & \cdots & \theta a_{1,M} + (1-\theta) m_{1,M} \\ \theta a_{2,1} + (1-\theta) m_{2,1} & \theta a_{2,2} + (1-\theta) m_{2,2} & \cdots & \theta a_{2,M} + (1-\theta) m_{2,M} \\ & & \vdots & \\ \theta a_{M,1} + (1-\theta) m_{M,1} & \theta a_{M,2} + (1-\theta) m_{M,2} & \cdots & \theta a_{M,M} + (1-\theta) m_{M,M} \end{bmatrix}.$$

The matrix $A_{ex} = (\theta a_{l,k} + (1-\theta) m_{l,k})$ has exactly similar structure, with the only difference being that now $l, k \in \mathcal{E}_x$. There are no type changes from one class to another, hence $A_{mx,ex} = ((m_{l,k}))$, with $l \in \mathcal{M}_x$ and $k \in \mathcal{E}_x$. Our focus is to investigate the evolution of the number of shares of $\mathcal{E}_x$ class particles when started with a particle of $\mathcal{M}_x$.

Note that this kind of branching processes can model various real-world applications, we encountered one such example while studying content propagation over social networks[3]. We briefly describe the model of [3] to motivate this study before we proceed further.

## 4.1 An example in viral marketing

Consider a viral marketing problem in an online social network (OSNs) with a large number of users. In viral marketing, the content providers (CPs)/advertisers create contents that are appealing to the users (e.g., giving offers, discounts). Depending upon composition/quality of the post, quantified by the control variable $\eta \in [0,1]$, users spread the information about the content by sharing it with their friends/connections in the OSN. The better the post quality $\eta$, the broader the reach. In particular, we describe the model for the post propagation of two competing CPs, namely CP-1 and CP-2, advertising their products/services. Denote by $\eta_i$ the post quality factor of Post-$P_i$ which corresponds to the CP-$i$ where $i = 1, 2$. In an OSN, posts/contents on TLs appear at various levels (see figure 2) based on their newness, for instance, News Feed on Facebook. This reverse chronological appearance of the posts, in any user's page, is called 'timeline' (TL) which is an essential component influencing the content propagation. (Each user represented by its TL). In what follows, a TL may hold: i) both Post-$P_1$ and Post-$P_2$, ii) Post-$P_1$ only, iii) Post-$P_2$ only, or iv) neither Post-$P_1$ nor Post-$P_2$ till, say, first $N$ levels. We study the time evolution of such TLs containing posts of interest. We use the following notations:

- $X_{l,k}(t)$ : denotes the number of TLs having Post-$P_1$ at level $l$ and Post-$P_2$ at level $k$ at time $t$

- $X_{l,0}(t)$ : denotes the number of TLs having Post-$P_1$ at level $l$, (and these TLs do not contain Post-$P_2$) at time $t$; and similarly $X_{0,k}(t)$ is defined.

The TLs are categorized in the following way:

- **Mixed class TLs**: the TLs holding both the Post-$P_1$ and Post-$P_2$ and referred as "mx" TLs.

- **Exclusive class-$i$ TLs**: the TLs holding Post-$i$ only and referred as "ex-$i$" TLs where $i = 1, 2$.

**Dynamics of content propagation**: A general mechanism for the propagation of the Post-P (say) is as follows. Primarily, a TL holding the Post-P either receives a post from the other TL (shift transition), or it shares the Post-P, possibly along with the other posts on it, with a random number of friends (share transition). In the former case, the position of the post of interest gets changed (shifts down by few lines), while the later case gives rise to new TLs with the post of interest. The schematic diagram 1 describes the transitions and the TL structure.

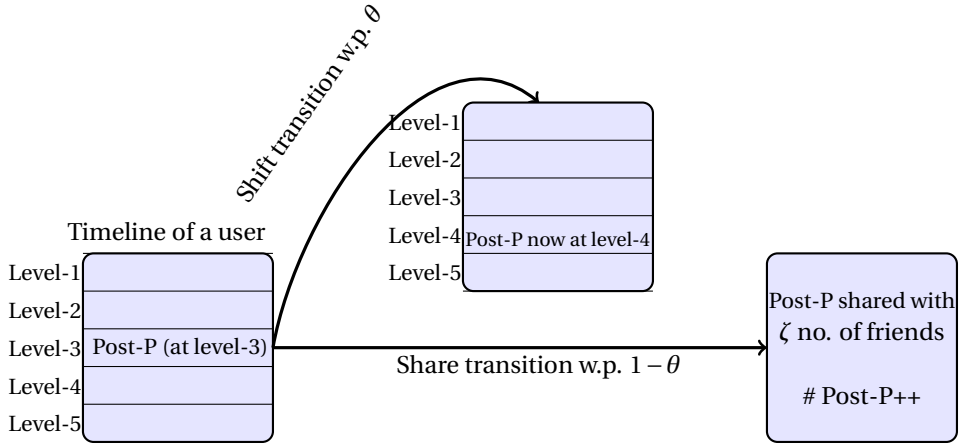

Figure 1: Propagation of Post-P: transitions

In the share transition, a TL holding Post-P, say at level $l$, wakes up after exponentially distributed time with parameter $\nu$, i.e., $exp(\nu)$. It reads the Post-P with probability $r_l$ and shares it with a random number of friends $\zeta$ depending Post-P quality factor $\eta$, and it produces type-$i$ ($i \le N$) with probability $\bar{\rho}_i$. Whereas in shift transition, one of the TLs holding Post-P receives a post after $exp(\lambda)$ time, which then pushes down the existing content on these TLs by level one. Observe that the probability of a shift is $\theta := \frac{\lambda}{\lambda+\nu}$ and share transition occurs with probability $1-\theta$. The diagram 2 describes the content propagation in the share transition. Note that the competing posts propagate through similar way, and we study the branching process modeling the competing posts.

**Branching model:** We model the above content propagation phenomenon as *continuous time multi-type branching process* (CMTBP). A TL with Post-$P_1$ (Post-$P_2$) can share it with random number of friends $\zeta_1/\zeta_2$, and produces only 'ex-1'('ex-2') types offspring respectively. On the other hand, a TL with both the posts can produce: i) 'ex-1' type when it shares Post-$P_1$ not Post-$P_2$ and vice-versa, ii) 'mx' type when it shares both Post-$P_1$ and Post-$P_2$. As the posts of competing CPs are similar, it is assumed the competing post lying at a lower level compared to the other receives reduced attention, and this effect is captured by the parameter $\delta \in [0,1]$.

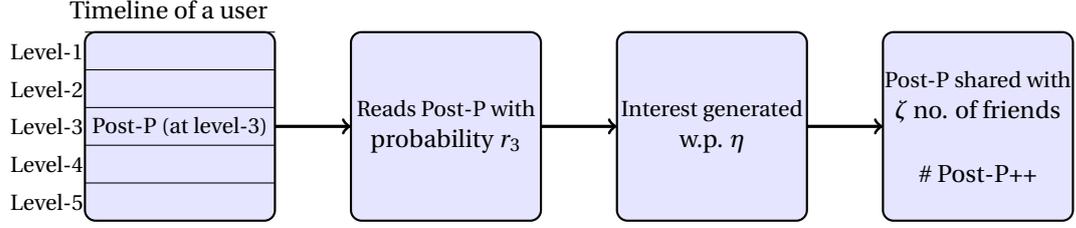

Figure 2: Share transition

Below we present the matrix $\mathbb{A}_{tc}$ (see [3] for more details) describing the transitions among these classes with $\mathbf{0}_j$ as the null matrix of appropriate order where $j = 1, \cdots, 4$, we have (see [3] for details)

$$\mathbb{A}_{tc} = (\lambda + \nu) \begin{bmatrix} \mathscr{A}_{mx} & A^1_{mx,ex} & A^2_{mx,ex} \\ \mathbf{0}_1 & A^1_{ex} & \mathbf{0}_2 \\ \mathbf{0}_3 & \mathbf{0}_4 & A^2_{ex} \end{bmatrix} ; \text{ where } A^i_{mx,ex} \text{ for } i=1,2, and$$

$$A^i_{mx,ex} = \begin{bmatrix} c_{mx,i} r_1 \bar{\rho}_1 & c_{mx,i} r_1 \bar{\rho}_2 & \cdots & c_{mx,i} r_1 \bar{\rho}_{N-1} & 0 \\ c'_{mx,i} r_1 \bar{\rho}_1 & c'_{mx,i} r_1 \bar{\rho}_2 & \cdots & c'_{mx,i} r_1 \bar{\rho}_{N-1} & 0 \\ c_{mx,i} r_2 \bar{\rho}_1 & c_{mx,i} r_2 \bar{\rho}_2 & \cdots & c_{mx,i} r_2 \bar{\rho}_{N-1} & 0 \\ c'_{mx,i} r_2 \bar{\rho}_1 & c'_{mx,i} r_2 \bar{\rho}_2 & \cdots & c'_{mx,i} r_2 \bar{\rho}_{N-1} & 0 \\ \vdots & \vdots & \ddots & \vdots & \vdots \\ c_{mx,i} r_{N-2} \bar{\rho}_1 & c_{mx,i} r_{N-2} \bar{\rho}_2 & \cdots & c_{mx,i} r_{N-2} \bar{\rho}_{N-2} & 0 \\ c'_{mx,i} r_{N-2} \bar{\rho}_1 & c'_{mx,i} r_{N-2} \bar{\rho}_2 & \cdots & c'_{mx,i} r_{N-1} \bar{\rho}_{N-2} & 0 \\ c_{mx,i} r_{N-1} \bar{\rho}_1 & c_{mx,i} r_{N-1} \bar{\rho}_2 & \cdots & c_{mx,i} r_{N-1} \bar{\rho}_{N-1} & \theta \\ c'_{mx,i} r_{N-1} \bar{\rho}_1 & c'_{mx,i} r_{N-1} \bar{\rho}_2 & \cdots & c'_{mx,i} r_{N-1} \bar{\rho}_{N-1} & 0 \end{bmatrix}, \mathscr{A}_{mx} = \begin{bmatrix} z'_1 r_1 - 1 & z_1 r_1 & \theta + z'_2 r_1 & \cdots & z'_{N-1} r_1 & z_{N-1} r_1 \\ z_1 r_1 & z'_1 r_1 - 1 & z_2 r_1 & \cdots & z_{N-1} r_1 & z'_{N-1} r_1 \\ z'_1 r_2 & z_1 r_2 & z'_2 r_2 - 1 & \cdots & z'_{N-1} r_2 & z_{N-1} r_2 \\ z_1 r_2 & z'_1 r_2 & z_2 r_2 & \cdots & z_{N-1} r_2 & z'_{N-1} r_2 \\ \vdots & \vdots & \vdots & \ddots & \vdots & \vdots \\ z'_1 r_{N-2} & z_1 r_{N-2} & z'_2 r_{N-2} & \cdots & \theta + z'_{N-1} r_{N-2} & z_{N-1} r_{N-2} \\ z_1 r_{N-2} & z'_1 r_{N-2} & z_2 r_{N-2} & \cdots & z_{N-1} r_{N-2} & \theta + z'_{N-1} r_{N-2} \\ z'_1 r_{N-1} & z_1 r_{N-1} & z'_1 r_{N-1} & \cdots & z'_{N-1} r_{N-1} - 1 & z_{N-1} r_{N-1} \\ z_1 r_{N-1} & z'_1 r_{N-1} & z_2 r_{N-1} & \cdots & z_{N-1} r_{N-1} & z'_{N-1} r_{N-1} - 1 \end{bmatrix}$$

where

- $c_{mx} = \delta(1-\theta)\eta_1\eta_2 m$, , $z'_j = (1-p)c_{mx}\bar{\rho}_j$ and $z_j = pc_{mx}\bar{\rho}_j$ for all $j$, and $m$ is the mean number of friends
- $c_{mx,i} = (1-\theta)m\eta_i(1-\delta\eta_{-i})$, $c'_{mx,i} = (1-\theta)m\delta\eta_i(1-\eta_{-i})$ $i=1,2$ and $-i = 1\mathbb{1}_{\{i=2\}} + 2\mathbb{1}_{\{i=1\}}$.

**Remarks:** We can observe two important aspects of this CMTBP:

- The *underlying CMTBP modeling content propagation of competing posts turns out to be decomposable*, a type of BP studied in previous sections;

- When a TL undergoes a shift transition (which happens with probability $\theta$), we have a different type of TL with the same post (now in different level); thus we have a type-changing process. This is indeed *type-changing vector decomposable branching process (TC-VDBP)*.

- We see that the analysis of VDCBP follows immediately as it is a specific case of TC-VDBP. To be more precise, $\theta = 0$ corresponds to this special case.

Observe that 'mx' class particles produce the offspring of all the classes whereas an 'ex' class particle produces offspring of its class only. This allows us to split the generator matrix into two sub-matrices as below, which would facilitate the required analysis:

$$\begin{bmatrix} \mathscr{A}_{mx} & A^1_{mx,ex} \\ \mathbf{0} & A^1_{ex} \end{bmatrix} \quad \text{and} \quad \begin{bmatrix} \mathscr{A}_{mx} & A^2_{mx,ex} \\ \mathbf{0} & A^2_{ex} \end{bmatrix}. \tag{14}$$

Note that evolution of particles corresponding to each sub-matrix can be studied separately and also observe that each of the sub-matrices corresponds to a TC-VDBP (given by (13)).

## 4.2 Analysis: time evolution of the expected number of shares

We are now ready to study the time evolution of the expected number of shares in TC-VDBP.

**Two different notions of total progeny:** We emphasize that there are two different notions for the total progeny in this peculiar TC-VDBP, as opposed to the standard one in the branching processes. One may view the type-changing as the production of one offspring of a different type, and thereby adding one to the total progeny (for each type-change). This phenomenon is the usual way the total progeny is counted in standard BPs. Alternatively, one may not view type-change as an offspring, which can lead to a different (new) notion of total progeny that counts only the new offspring. For instance, as we already discussed, in a social network one needs only the count of total shares (the number of distinct users shared with the post of interest). *Inspired by social network terminology, we refer to the progeny that does not count the type-changes as 'number of shares,' while the one that counts all the transitions as the usual 'total progeny.'*

In this section, we derive the time evolution of the expected number of shares. It is clear that the total progeny of VDCBP is obtained by substituting $\theta = 0$ in the expression for the expected number of shares.

The number of shares at a time instance, say $t$, represents the total accumulated population (i.e., including the dying particles) of all types till $t$. We study the evolution of the number of shares for both exclusive class and mixed class particles. Note that the evolution of the population of exclusive/mixed class when initiated with its own class particle(s) is obtained using the well-known theory of non-decomposable BPs (e.g., see [4]). The specific expression for the number of shares can be found in [3]. According to which the number of shares say $y_l^e(t)$ of an $\mathcal{E}_x$ class till time $t$ when initiated with its own class type-$l$ particle is: (see [3] for more details)

$$y_l^e(t) = \tilde{g}_l^e + h_l^e e^{\alpha_e t}, \quad \text{with} \quad \tilde{g}_l^e = 1 + (\lambda + \nu)\left\{A_{ex}^{-1}\mathbf{k}\right\}_l, \quad h_l^e = -(\lambda + \nu)\left\{A_{ex}^{-1}\mathbf{k}\right\}_l; \quad where$$

- $\left\{A_{ex}^{-1}\mathbf{k}\right\}_l$ denote the $l-th$ component of the vector $A_{ex}^{-1}\mathbf{k}$;

- $\mathbf{k} = [1-\theta, 1-\theta, \cdots, 1-\theta, 1]^T$, and $\alpha_e$ is largest the eigenvalue of the matrix $A_{ex}$.

Note that we adopted the convention that $y_l^e(0) = 1$ for $l \in \mathcal{E}_x$. For the sake of convenience, one can take $y_l^e(0) = 0$ to study the problem without loss of generality. Basically, depending on the context of the problem, one may or may not count the starting particle as one share added. In what follows, we now have

$$y_l^e(t) = g_l^e + h_l^e e^{\alpha_e t}, \quad \text{with} \quad g_l^e = -h_l^e \quad \text{for } l \in \mathcal{E}_x, \tag{15}$$

We now focus on investigating the evolution of the expected number of shares of each exclusive class when the process starts with a mixed class particle. We obtain this by first deriving appropriate fixed point equations.

### 4.2.1 Derivation of an appropriate fixed point (FP) equation

Denote by $y_l(t)$ the number of shares of $\mathcal{E}_x$ class till time $t$ when the process is initiated with a type-$l$ particle of the $\mathcal{M}_x$ class, and $\mathbf{y}(t) = \{y_l(t)\}_l$ represents the 'expected shares' vector. We conjecture that $y_l(t)$ satisfies a fixed point equation in an appropriate functional space, i.e., $z_l(t) = G(\mathbf{z}(t))$ where $\mathbf{z}(t) := \{z_l(t)\}_l$ represent finite number of waveforms on time interval $[0, \infty)$. We arrive at the fixed point equation by conditioning on the events related to the first transition epoch. Let the random variable $\tau$ represent the time instance of the first transition epoch, which is exponentially distributed with parameter $\lambda + \nu$. Conditioning on the first transition events, we observe that the number of shares $\mathbf{y}(\cdot)$ satisfy the following fixed point equation:

$$\begin{aligned} y_l(t) &= \theta \int_0^t \sum_{k \in \mathcal{M}_x} a_{l,k} y_k(t-\tau)(\lambda+\nu)e^{-(\lambda+\nu)\tau} d\tau \\ &+ (1-\theta)\int_0^t \sum_{k \in \mathcal{M}_x} m_{l,k}(1+y_k(t-\tau))(\lambda+\nu)e^{-(\lambda+\nu)\tau}d\tau \\ &+ (1-\theta)\int_0^t \sum_{k \in \mathcal{E}_x} m_{l,k}(1+y_k^e(t-\tau))(\lambda+\nu)e^{-(\lambda+\nu)\tau}d\tau. \end{aligned}$$

The above is due to the following reasons:

- The type-$l$ undergoes a shift transition w.p. $\theta$, its type gets changed to type-$k$ of the same class (i.e., $\mathcal{M}_x$).

- The type-$l$ undergoes a share transition w.p. $1-\theta$ it produces $m_{l,k}$ offspring belonging to either class $\mathcal{M}_x$ or $\mathcal{E}_x$. As per example, it produces particles of $\mathcal{M}_x$ when both Post-$P_1$ and Post-$P_2$ are shared, whereas particles of $\mathcal{E}_x$ are produced when only one the posts is shared.

### 4.2.2 Solution of the fixed point equation

We assume the following structure for fixed point waveform, $y_l(t) = g_l + h_l e^{\alpha_e t} + o_l e^{\bar{\alpha} t}$ for $l \in \mathcal{M}_x$ and $\bar{\alpha}$ is a constant (which we will find out). We show that these kind of functions indeed satisfy the appropriate fixed point equations. However, we are yet to investigate that the fixed point solution of the $z_l(t) = G(\mathbf{z}(t))$ is unique, and hence, equals $y_l(t)$ as suggested above. We now derive a solution of the above fixed point equation. Towards this, we have the following Lemma. Let $\alpha_e, \alpha_{mx}$ be the largest eigenvalue of the matrices $(A_{ex} - I)\lambda_v, (A_{mx} - I)\lambda_v$ respectively.

**Theorem 4** *When $\alpha_e > 0$, i.e. the exclusive class is super-critical, a solution of the above fixed point equation $y_l(t) = G_l(\mathbf{y}(t))$ is the following:*

1. *When the $\mathcal{M}_x$ population gets extinct with probability one (i.e., when $\alpha_{mx} < 0$), then $y_l(t) = y_l^e(t) = g_l^e + h_l^e e^{\alpha_e t}$, with $g_l^e = -h_l^e$.*

2. *When the $\mathcal{M}_x$ population survives with non zero probability (i.e., when $\alpha_{mx} > 0$), then*

$$y_l(t) = g_l + h_l e^{\alpha_e t} + o_l e^{\bar{\alpha} t} \qquad (16)$$

*where $g_l, h_l, o_l$ are as given as:*

$$
\begin{aligned}
h_l &= \frac{\lambda_v(1-\theta)\sum_{k\in\mathcal{E}_x\cup\mathcal{M}_x} m_{l,k}}{\alpha_e} - \frac{(1-\theta)\lambda_v \sum_{k\in\mathcal{E}_x} m_{l,k} h_k^e}{(\bar{\alpha}-\alpha_e)} + \frac{(1-\theta)\lambda_v}{(\bar{\alpha}-\alpha_e)} \sum_{k\in\mathcal{M}_x\cup\mathcal{E}_x} \bar{\alpha}(\lambda_v+\alpha_e)m_{l,k} \\
&\quad - \frac{(1-\theta)\lambda_v\bar{\alpha}}{(\bar{\alpha}-\alpha_e)} \left( \lambda_v \sum_{k\in\mathcal{M}_x} \left(\theta a_{l,k}+(1-\theta)m_{l,k}\right) \sum_{k'\in\mathcal{E}_x\cup\mathcal{M}_x} m_{k,k'} \right)
\end{aligned} \qquad (17)
$$

$$
\begin{aligned}
g_l &= -\frac{\lambda_v(1-\theta)\sum_{k\in\mathcal{E}_x\cup\mathcal{M}_x} m_{l,k}}{\alpha_e} + \frac{(1-\theta)\lambda_v\sum_{k\in\mathcal{E}_x} m_{l,k}h_k^e}{\bar{\alpha}} - \frac{(1-\theta)\lambda_v}{\bar{\alpha}\alpha_e}\sum_{k\in\mathcal{M}_x\cup\mathcal{E}_x}(\lambda_v+\alpha_e)m_{l,k} \\
&\quad + \frac{(1-\theta)\lambda_v}{\bar{\alpha}\alpha_e}\left(\lambda_v\sum_{k\in\mathcal{M}_x}\left(\theta a_{l,k}+(1-\theta)m_{l,k}\right)\sum_{k'\in\mathcal{E}_x\cup\mathcal{M}_x} m_{k,k'}\right).
\end{aligned} \qquad (18)
$$

$$
\begin{aligned}
o_l &= \frac{(1-\theta)(\lambda_v+\alpha_e)\lambda_v}{(\bar{\alpha}-\alpha_e)\bar{\alpha}}\left( \frac{\alpha_e\sum_{k\in\mathcal{E}_x} m_{l,k}h_k^e}{\lambda_v+\alpha_e} - \sum_{k\in\mathcal{M}_x\cup\mathcal{E}_x} m_{l,k} + \frac{\lambda_v\sum_{k\in\mathcal{M}_x}\left(\theta a_{l,k}+(1-\theta)m_{l,k}\right)\sum_{k'\in\mathcal{E}_x\cup\mathcal{M}_x}m_{k,k'}}{\lambda_v+\alpha_e} \right) \\
&= \frac{(1-\theta)\lambda_v\alpha_e\sum_{k\in\mathcal{E}_x}m_{l,k}h_k^e}{(\bar{\alpha}-\alpha_e)\bar{\alpha}} + \frac{(1-\theta)\lambda_v}{(\bar{\alpha}-\alpha_e)\bar{\alpha}}\left(\lambda_v\sum_{k\in\mathcal{M}_x}\left(\theta a_{l,k}+(1-\theta)m_{l,k}\right)\sum_{k'\in\mathcal{E}_x\cup\mathcal{M}_x}m_{k,k'} - (\lambda_v+\alpha_e)\sum_{k\in\mathcal{M}_x\cup\mathcal{E}_x}m_{l,k}\right)
\end{aligned} \qquad (19)
$$

*and $\bar{\alpha}$ is as given in equation (24), i.e., $\bar{\alpha} = (eig(A_{mx})-1)\lambda_v$. Assume that $\alpha_{mx}$ is the only eigenvalue of $(A_{mx}-I)\lambda_v$ larger than zero, then we have $\bar{\alpha} = \alpha_{mx}$.*

**Proof 4** *As $y_l(t)$ satisfies the following fixed point equation*

$$
\begin{aligned}
y_l(t) &= \theta\int_0^t \sum_{k\in\mathcal{M}_x} a_{l,k} y_k(t-\tau)(\lambda+\nu)e^{-(\lambda+\nu)\tau} d\tau \\
&\quad +(1-\theta)\int_0^t \sum_{k\in\mathcal{M}_x} m_{l,k}(1+y_k(t-\tau))(\lambda+\nu)e^{-(\lambda+\nu)\tau} d\tau \\
&\quad +(1-\theta)\int_0^t \sum_{k\in\mathcal{E}_x} m_{l,k}(1+y_k^e(t-\tau))(\lambda+\nu)e^{-(\lambda+\nu)\tau} d\tau
\end{aligned}
$$

*and that it has the following structure: $y_l(t) = g_l + h_l e^{\alpha_e t} + o_l e^{\bar{\alpha} t}$ for $l\in\mathcal{M}_x$.*

*Directly substituting the above representation of $y_l(t)$, we have ($\lambda_\nu := \lambda + \nu$):*

$$
\begin{aligned}
g_l + h_l e^{\alpha_e t} + o_l e^{\bar{\alpha} t} &= y_l(t) = \theta \int_0^t \sum_{k \in \mathcal{M}_x} a_{l,k} \left( g_k + h_k e^{\alpha_e(t-\tau)} + o_k e^{\bar{\alpha}(t-\tau)} \right) \lambda_\nu e^{-\lambda_\nu \tau} d\tau \\
&\quad + (1-\theta) \int_0^t \sum_{k \in \mathcal{M}_x} m_{l,k} \left( 1 + g_k + h_k e^{\alpha_e(t-\tau)} + o_k e^{\bar{\alpha}(t-\tau)} \right) \lambda_\nu e^{-\lambda_\nu \tau} d\tau \\
&\quad + (1-\theta) \int_0^t \sum_{k \in \mathcal{E}_x} m_{l,k} \left( 1 + \left( g_k^e + h_k^e e^{\alpha_e(t-\tau)} \right) \right) \lambda_\nu e^{-\lambda_\nu \tau} d\tau \\
&= \left[ \theta \sum_{k \in \mathcal{M}_x} a_{l,k} g_k + (1-\theta) \sum_{k \in \mathcal{M}_x} m_{l,k}(1 + g_k) + (1-\theta) \sum_{k \in \mathcal{E}_x} m_{l,k}(1 + g_k^e) \right] (1 - e^{-\lambda_\nu t}) \\
&\quad + e^{\alpha_e t} \left[ \theta \sum_{k \in \mathcal{M}_x} a_{l,k} h_k + (1-\theta) \sum_{k \in \mathcal{M}_x} m_{l,k} h_k + (1-\theta) \sum_{k \in \mathcal{E}_x} m_{l,k} h_k^e \right] \left( 1 - e^{-(\lambda_\nu + \alpha_e) t} \right) \frac{\lambda_\nu}{\lambda_\nu + \alpha_e} \\
&\quad + e^{\bar{\alpha} t} \left[ \theta \sum_{k \in \mathcal{M}_x} a_{l,k} o_k + (1-\theta) \sum_{k \in \mathcal{M}_x} m_{l,k} o_k \right] \left( 1 - e^{-(\lambda_\nu + \bar{\alpha}) t} \right) \frac{\lambda_\nu}{\lambda_\nu + \bar{\alpha}}.
\end{aligned}
$$

*Thus we need that:*

$$
g_l = \theta \sum_{k \in \mathcal{M}_x} a_{l,k} g_k + (1-\theta) \sum_{k \in \mathcal{M}_x} m_{l,k}(1 + g_k) + (1-\theta) \sum_{k \in \mathcal{E}_x} m_{l,k}(1 + g_k^e), \tag{20}
$$

$$
h_l = \left( \theta \sum_{k \in \mathcal{M}_x} a_{l,k} h_k + (1-\theta) \sum_{k \in \mathcal{M}_x} m_{l,k} h_k + (1-\theta) \sum_{k \in \mathcal{E}_x} m_{l,k} h_k^e \right) \frac{\lambda_\nu}{\lambda_\nu + \alpha_e} \quad \text{and} \tag{21}
$$

$$
o_l = \left( \theta \sum_{k \in \mathcal{M}_x} a_{l,k} o_k + (1-\theta) \sum_{k \in \mathcal{M}_x} m_{l,k} o_k \right) \frac{\lambda_\nu}{\lambda_\nu + \bar{\alpha}} \tag{22}
$$

*and that*

$$
\begin{aligned}
&-\theta \sum_{k \in \mathcal{M}_x} a_{l,k} g_k - (1-\theta) \sum_{k \in \mathcal{M}_x} m_{l,k}(1 + g_k) - (1-\theta) \sum_{k \in \mathcal{E}_x} m_{l,k}(1 + g_k^e) \\
&- \left[ \theta \sum_{k \in \mathcal{M}_x} a_{l,k} h_k + (1-\theta) \sum_{k \in \mathcal{M}_x} m_{l,k} h_k + (1-\theta) \sum_{k \in \mathcal{E}_x} m_{l,k} h_k^e \right] \frac{\lambda_\nu}{\alpha_e + \lambda_\nu} \\
&- \left[ \theta \sum_{k \in \mathcal{M}_x} a_{l,k} o_k + (1-\theta) \sum_{k \in \mathcal{M}_x} m_{l,k} o_k \right] \frac{\lambda_\nu}{\bar{\alpha} + \lambda_\nu} = 0
\end{aligned}
$$

*which directly follows as we have $y(0) = 0$, i.e., the user with post of interest has not shared, and thus*

$$
g_l + h_l + o_l = 0.
$$

**Case I:** *When mixed gets extinct with probability one, there is no positive $\bar{\alpha}$, that satisfies the eigen value equation of hypothesis, also given by equation (24). In this case, we will show that $o_l = 0$ and $g_l = -h_l$ is the solution. Note that in this case $h_l = h_l^e$ (see equation 15). Now using (20) and (21), we have:*

$$
\begin{aligned}
g_l + \frac{\alpha_e + \lambda_\nu}{\lambda_\nu} h_l &= \theta \sum_{k \in \mathcal{M}_x} a_{l,k}(g_k + h_k) + (1-\theta) \sum_{k \in \mathcal{M}_x} m_{l,k}(1 + g_k + h_k) + (1-\theta) \sum_{k \in \mathcal{E}_x} m_{l,k}(1 + g_k^e + h_k^e) \\
&= (1-\theta) \sum_{k \in \mathcal{M}_x} m_{l,k} + (1-\theta) \sum_{k \in \mathcal{E}_x} m_{l,k}, \quad \text{and} \\
h_l &= \left( -g_l + (1-\theta) \sum_{k \in \mathcal{E}_x \cup \mathcal{M}_x} m_{l,k} \right) \frac{\lambda_\nu}{\lambda_\nu + \alpha_e} = \frac{\lambda_\nu (1-\theta) \sum_{k \in \mathcal{E}_x \cup \mathcal{M}_x} m_{l,k}}{\alpha_e}.
\end{aligned}
$$

*Thus,* $y_l(t) = y_l^e(t) = -\dfrac{\lambda_\nu (1-\theta) \sum_{k \in \mathcal{E}_x \cup \mathcal{M}_x} m_{l,k}}{\alpha_e} + \dfrac{\lambda_\nu (1-\theta) \sum_{k \in \mathcal{E}_x \cup \mathcal{M}_x} m_{l,k}}{\alpha_e} e^{\alpha_e t}$

**Case II:** *When $o_l \neq 0$, then using $g_l + h_l + o_l = 0$, $g_l^e + h_l^e = 0$ for each $l$ (note that here $h_l \neq h_l^e$ ). And using (20)-(22) we get*

$$g_l + \frac{\lambda_v + \alpha_e}{\lambda_v} h_l + \frac{\lambda_v + \bar{\alpha}}{\lambda_v} o_l = \theta \sum_{k \in \mathcal{M}_x} a_{l,k}(g_k + h_k + o_k) + \sum_{k \in \mathcal{M}_x} (1-\theta) m_{l,k}(1 + g_k + h_k + o_k)$$
$$+ \sum_{k \in \mathcal{E}_x} (1-\theta) m_{l,k}\left(1 + g_k^e + h_k^e\right) = \sum_{k \in \mathcal{E}_x \cup \mathcal{M}_x} (1-\theta) m_{l,k}$$
$$h_l = \left[-g_l + (1-\theta) \sum_{k \in \mathcal{E}_x \cup \mathcal{M}_x} m_{l,k} - o_l \frac{\bar{\alpha} + \lambda_v}{\lambda_v}\right] \frac{\lambda_v}{\lambda_v + \alpha_e}$$
$$= \left[-g_l - o_l + (1-\theta) \sum_{k \in \mathcal{E}_x \cup \mathcal{M}_x} m_{l,k} - o_l \frac{\bar{\alpha}}{\lambda_v}\right] \frac{\lambda_v}{\lambda_v + \alpha_e}.$$

Thus

$$h_l = \frac{\lambda_v(1-\theta) \sum_{k \in \mathcal{E}_x \cup \mathcal{M}_x} m_{l,k} - o_l \bar{\alpha}}{\alpha_e}, \text{ and } g_l = -h_l - o_l = -\frac{\lambda_v(1-\theta) \sum_{k \in \mathcal{E}_x \cup \mathcal{M}_x} m_{l,k} - o_l(\bar{\alpha} - \alpha_e)}{\alpha_e}. \quad (23)$$

Now summing equations (20)-(22), and by using $h_l + g_l + o_l = 0$ and $h_l^e = -g_l^e$ in (21), we have:

$$0 = g_l + h_l + o_l = \left[\theta \sum_{k \in \mathcal{M}_x} a_{l,k} g_k + (1-\theta) \sum_{k \in \mathcal{M}_x} m_{l,k}(1 + g_k) + (1-\theta) \sum_{k \in \mathcal{E}_x} m_{l,k}(1 + g_k^e)\right]$$
$$+ \left[\theta \sum_{k \in \mathcal{M}_x} a_{l,k} h_k + (1-\theta) \sum_{k \in \mathcal{M}_x} m_{l,k} h_k + (1-\theta) \sum_{k \in \mathcal{E}_x} m_{l,k} h_k^e\right]\left(1 - \frac{\alpha_e}{\lambda_v + \alpha_e}\right)$$
$$+ \left[\theta \sum_{k \in \mathcal{M}_x} a_{l,k} o_k + (1-\theta) \sum_{k \in \mathcal{M}_x} m_{l,k} o_k\right]\left(1 - \frac{\bar{\alpha}}{\lambda_v + \bar{\alpha}}\right)$$
$$= -\left[\theta \sum_{k \in \mathcal{M}_x} a_{l,k} h_k + (1-\theta) \sum_{k \in \mathcal{M}_x} m_{l,k} h_k + (1-\theta) \sum_{k \in \mathcal{E}_x} m_{l,k} h_k^e\right] \frac{\alpha_e}{\lambda_v + \alpha_e}$$
$$- \left[\theta \sum_{k \in \mathcal{M}_x} a_{l,k} o_k + (1-\theta) \sum_{k \in \mathcal{M}_x} m_{l,k} o_k\right] \frac{\bar{\alpha}}{\lambda_v + \bar{\alpha}} + (1-\theta) \sum_{k \in \mathcal{M}_x \cup \mathcal{E}_x} m_{l,k}.$$

Substituting the value of $h_k$ as given in equation (23), we get:

$$-\sum_{k \in \mathcal{M}_x} \left(\theta a_{l,k} + (1-\theta) m_{l,k}\right) \left(\frac{\lambda_v(1-\theta) \sum_{k' \in \mathcal{E}_x \cup \mathcal{M}_x} m_{l,k'} - o_l \bar{\alpha}}{\alpha_e}\right) \times \frac{\alpha_e}{\lambda_v + \alpha_e} - (1-\theta) \sum_{k \in \mathcal{E}_x} m_{l,k} h_k^e \frac{\alpha_e}{\lambda_v + \alpha_e}$$
$$- \left[\theta \sum_{k \in \mathcal{M}_x} a_{l,k} o_k + (1-\theta) \sum_{k \in \mathcal{M}_x} m_{l,k} o_k\right] \frac{\bar{\alpha}}{\lambda_v + \bar{\alpha}} + (1-\theta) \sum_{k \in \mathcal{M}_x \cup \mathcal{E}_x} m_{l,k} = 0$$

$$\left[\theta \sum_{k \in \mathcal{M}_x} a_{l,k} o_k + (1-\theta) \sum_{k \in \mathcal{M}_x} m_{l,k} o_k\right]\left(\frac{\bar{\alpha}}{\lambda_v + \alpha_e} - \frac{\bar{\alpha}}{\lambda_v + \bar{\alpha}}\right) = (1-\theta) \sum_{k \in \mathcal{E}_x} m_{l,k} h_k^e \frac{\alpha_e}{\lambda_v + \alpha_e} - (1-\theta) \sum_{k \in \mathcal{M}_x \cup \mathcal{E}_x} m_{l,k}$$
$$+ \left[\sum_{k \in \mathcal{M}_x} \left(\theta a_{l,k} + (1-\theta) m_{l,k}\right) \sum_{k' \in \mathcal{E}_x \cup \mathcal{M}_x} (1-\theta) m_{k,k'}\right] \frac{\lambda_v}{\lambda_v + \alpha_e}$$
$$o_l \frac{\lambda_v + \bar{\alpha}}{\lambda_v}\left(\frac{\bar{\alpha}}{\lambda_v + \alpha_e} - \frac{\bar{\alpha}}{\lambda_v + \bar{\alpha}}\right) = 1-\theta) \sum_{k \in \mathcal{E}_x} m_{l,k} h_k^e \frac{\alpha_e}{\lambda_v + \alpha_e} - (1-\theta) \sum_{k \in \mathcal{M}_x \cup \mathcal{E}_x} m_{l,k}$$
$$+ (1-\theta) \left[\sum_{k \in \mathcal{M}_x} \left(\theta a_{l,k} + (1-\theta) m_{l,k}\right) \sum_{k' \in \mathcal{E}_x \cup \mathcal{M}_x} m_{k,k'}\right] \frac{\lambda_v}{\lambda_v + \alpha_e} \quad \text{(using equation (22)).}$$

Thus $o_l$ is given by (19) given in the hypothesis of the Lemma.

Similarly, one can obtain the close form expressions for $h_l$ and $g_l$ by substituting the value of $o_l$ in (23) and the final expressions are as in (17)-(18) given in the hypothesis of the Lemma.

*Again from (22), observe that the vector $\mathbf{o}_l := \{o_l\}_{l \in \mathcal{M}_x}$ is an eigenvector of the following matrix corresponding to eigenvalue one:*

$$\frac{\lambda_v}{\lambda_v + \bar{\alpha}} \mathcal{A}_o, \text{ with } \mathcal{A}_o := \begin{bmatrix} \theta a_{1,1} + (1-\theta) m_{1,1} & \theta a_{1,2} + (1-\theta) m_{1,2} & \cdots & \theta a_{1,M} + (1-\theta) m_{1,M} \\ \theta a_{2,1} + (1-\theta) m_{2,1} & \theta a_{2,2} + (1-\theta) m_{2,2} & \cdots & \theta a_{2,M} + (1-\theta) m_{2,M} \\ & & \vdots & \\ \theta a_{M,1} + (1-\theta) m_{M,1} & \theta a_{M,2} + (1-\theta) m_{M,2} & \cdots & \theta a_{M,M} + (1-\theta) m_{M,M} \end{bmatrix},$$

*which exactly equals $A_{mx}$. In other words, we have $\frac{\lambda_v}{\lambda_v + \bar{\alpha}} \mathcal{A}_o \mathbf{o}_l = 1 \times \mathbf{o}_l$. Thus, one must chose a positive $\bar{\alpha}$ such that the eigen value of $\frac{\lambda_v}{\lambda_v + \bar{\alpha}} \mathcal{A}_o$ is one, or in other words:*

$$\bar{\alpha} = (eig(\mathcal{A}_o) - 1) \lambda_v. \tag{24}$$

*Also note that when $\bar{\alpha} > 0$, then $\mathbf{o}_l \neq \mathbf{0}$ where $\mathbf{0}$ is the zero vector of appropriate order. Thus, we have*

$$y_l(t) = g_l + h_l e^{\alpha_e t} + o_l e^{\bar{\alpha} t}; \text{ where } g_l, h_l, o_l \text{ are as given in by the equations 18, 17, 19.}$$

∎

**Remarks:** 1) When there are two or more eigenvalues of $(A_{mx} - I) \lambda_v$ larger than 0, then one may have more choices for $\bar{\alpha}$ in the above result. We need to work further to understand this type of scenarios. However, in the branching process literature (in other contexts) it is well-known that the largest eigenvalue determines the growth rate, and we anticipate the same here (in which case we will have $\bar{\alpha} = \alpha_{mx}$).

2) Thus, we notice that for decomposable branching processes, the grow rates of the current population (10) as well as the total shares (16) are influenced by two distinct exponential functions. Also, both of them (like in other variants of the branching process) are influenced by the same growth patterns. Total shares/total progeny even, in this case, does not seem to have a growth rate different from that of the current population.

## Conclusions

We study various variants of decomposable branching processes in the continuous time framework. In literature, decomposable branching processes are relatively less studied objects compared to the non-decomposable or irreducible counterpart. These objects are used in studying many real-world applications such as cancer biology, viral marketing, etc. We investigated the time evolution of the population of various classes when the process is initiated by the other class particle(s). We obtained various measures such as class-wise extinction probability, and the growth of the population of different classes in both scalar and vector versions of the decomposable BPs. We then studied another peculiar type of decomposable branching process where any parent at the transition epoch either produces a random number of offspring, or its type gets changed which may or may not be regarded as new offspring produced depending on the application. These processes arise in the modelling of content propagation of competing contents in online social networks. We obtained extinction probability, the growth rate of class-wise population for this case also. Additionally, we conjecture the time evolution of the expected number of shares (different from the total progeny in irreducible branching processes) as sum of two (distinct) exponential curves corresponding to the two classes.

# 5  Appendix

**Proof 5  of Theorem 1: (i) (a)**

With $\mathbb{J}_k = \{j_1, \cdots, j_k\}$ and $l < j_1 < \cdots < j_k < m$, we define the following set of indices:

$$S_l^m(\mathbb{J}_k) := \begin{cases} \{l, j_1, j_2, \cdots, j_k, m\} & \text{if } \mathbb{J}_k \neq \emptyset \\ \{l, m\} & \text{otherwise.} \end{cases}$$

For the ease of reference, we reproduce $M_m(t)$ given by (5) here and also write down the corresponding coefficients using the above definition:

$$\begin{aligned} M_m(t) &= X_m(t)e^{-\alpha_m t} + \sum_{k=1}^{m-1}\sum_{j_1=1}^{m-k}\sum_{j_2=j_1+1}^{m-k+1}\sum_{j_3=j_2+1}^{m-k+2}\cdots\sum_{j_k=j_{k-1}+1}^{m-1} \frac{\alpha_{j_1 j_2}\alpha_{j_2 j_3}\cdots\alpha_{j_k m} X_{j_1}(t)e^{-\alpha_m t}}{(\alpha_m - \alpha_{j_1})(\alpha_m - \alpha_{j_2})\cdots(\alpha_m - \alpha_{j_k})} \\ &= \left(a_1^m X_1(t) + a_2^m X_2(t) + \cdots + a_{m-1}^m X_{m-1}(t) + a_m^m X_m(t)\right)e^{-\alpha_m t}, \end{aligned} \quad (25)$$

where $a_i^i = 1\ \forall i$ and $a_i^m$ is the coefficient corresponding to term $X_i(t)e^{-\alpha_m t}$ in $M_m(t)$ and is given as:

$$\begin{aligned} a_i^m &= \sum_{k=2}^{m-i}\sum_{j_2=i+1}^{m-k+1}\cdots\sum_{j_k=j_{k-1}+1}^{m-1} \frac{\alpha_{i j_2}\alpha_{j_2 j_3}\cdots\alpha_{j_k m}}{(\alpha_m - \alpha_i)(\alpha_m - \alpha_{j_2})\cdots(\alpha_m - \alpha_{j_k})} + \frac{\alpha_{im}}{\alpha_m - \alpha_i}; \ i = 1, 2, \cdots, m-1 \\ &= \sum_{k=1}^{m-i-1}\sum_{j_1=i+1}^{m-k}\cdots\sum_{j_k=j_{k-1}+1}^{m-1} \frac{\alpha_{i j_1}\alpha_{j_1 j_2}\cdots\alpha_{j_k m}}{(\alpha_m - \alpha_i)(\alpha_m - \alpha_{j_2})\cdots(\alpha_m - \alpha_{j_k})} + \frac{\alpha_{im}}{\alpha_m - \alpha_i}. \end{aligned} \quad (26)$$

We need to prove the following: for any $\delta > 0$

$$E\left[M_m(t+\delta)\Big|\mathscr{F}_t\right] = E\left[\sum_{i=1}^{m} a_i^m X_i(t+\delta)e^{-\alpha_m(t+\delta)}\Big|\mathscr{F}_t\right] = M_m(t)\ \forall t;\ \text{recall } \mathscr{F}_t \text{ is the natural filtration.}$$

*Towards this, we first compute* $E[X_i(t+\delta)|\mathscr{F}_t]$. *Denote by* $(e^{B_m\delta})_{ij}$ *the* $(i,j)$-*th entry of the generator matrix*[5] $e^{B_m\delta}$ *(see equation (2)). Appealing to the equations (4), (5) and (6) as given in [4], we can write the following*

$$E\left[X_i(t+\delta)\big|\mathscr{F}_t\right] = \left(e^{B\delta}\right)_{1i}X_1(t) + \left(e^{B\delta}\right)_{2i}X_2(t) + \cdots + \left(e^{B\delta}\right)_{i-1\,i}X_{i-1}(t) + \left(e^{B\delta}\right)_{ii}X_i(t);\ \text{and}$$

$$E\left[M_m(t+\delta)\big|\mathscr{F}_t\right] = e^{-\alpha_m(t+\delta)}E\left[a_1^m X_1(t+\delta) + a_2^m X_2(t+\delta) + \cdots + a_m^m X_m(t+\delta)\big|\mathscr{F}_t\right]$$

$$= e^{-\alpha_m(t+\delta)}\left(\sum_{i=1}^{m}\left(e^{B\delta}\right)_{1i}a_i^m X_1(t) + \sum_{i=2}^{m}\left(e^{B\delta}\right)_{2i}a_i^m X_2(t) + \cdots +\right.$$

$$\left. + \sum_{i=m-1}^{m}\left(e^{B\delta}\right)_{m-1\,i}a_i^m X_{m-1}(t) + \left(e^{B\delta}\right)_{mm}a_m^m X_m(t)\right). \tag{27}$$

*Observe that if we prove* $\sum_{i=j}^{m}(e^{B_m\delta})_{ji}a_i^m = a_j^m e^{\alpha_m \delta}\ \forall\ j$, *then using equations (27) and (25), it follows that*
$E[M_m(t+\delta)|\mathscr{F}_t] = M_m(t)\ \forall\ \delta > 0$ *and hence the proof would be complete.*

**Claim**: *We have*[6] $\sum_{i=j}^{m}(e^{B_m\delta})_{ji}a_i^m = a_j^m e^{\alpha_m\delta}$ *for any* $j$.
*We begin with expanding the summation as follows:*

$$\sum_{i=j}^{m}\left(e^{B_m\delta}\right)_{ji}a_i^m = \left(e^{B_m\delta}\right)_{jj}a_j^m$$
$$+ \left(e^{B_m\delta}\right)_{j\,j+1}a_{j+1}$$
$$+ \cdots + \left(e^{B\delta}\right)_{jm}a_m^m, \tag{28}$$

*Using Lemma 1 and equation (26), these terms can be rewritten as*

$$\sum_{i=j}^{m}\left(e^{B_m\delta}\right)_{ji}a_i^m = e^{\alpha_j\delta}\left(\sum_{k=1}^{m-j-1}\sum_{j_1=j+1}^{m-k}\cdots\sum_{j_k=j_{k-1}+1}^{m-1}\frac{\alpha_{jj_1}\alpha_{j_1 j_2}\cdots\alpha_{j_k m}}{(\alpha_m-\alpha_j)(\alpha_m-\alpha_{j_1})\cdots(\alpha_m-\alpha_{j_k})} + \frac{\alpha_{jm}}{\alpha_m-\alpha_j}\right)$$

$$+ \alpha_{jj+1}\frac{e^{\alpha_j\delta}-e^{\alpha_{j+1}\delta}}{\alpha_j-\alpha_{j+1}}\left(\sum_{k=1}^{m-j-2}\sum_{j_1=j+2}^{m-k}\cdots\sum_{j_k=j_{k-1}+1}^{m-1}\frac{\alpha_{j+1\,j_1}\alpha_{j_1 j_2}\cdots\alpha_{j_k m}}{(\alpha_m-\alpha_{j+1})(\alpha_m-\alpha_{j_1})\cdots(\alpha_m-\alpha_{j_k})} + \frac{\alpha_{j+1\,m}}{\alpha_m-\alpha_{j+1}}\right)$$

$$+ \left(\alpha_{jj+2}\frac{e^{\alpha_j\delta}-e^{\alpha_{j+2}\delta}}{\alpha_j-\alpha_{j+2}} + \sum_{k=1}^{1}\sum_{j_1=j+1}^{j+1}\sum_{p\in S_j^{j+2}(\mathbb{J}_k)}\frac{\alpha_{jj_1}\alpha_{j_1 j+2}e^{\alpha_p\delta}}{\prod_{l\in S_j^{j+2}(\mathbb{J}_k):l\neq p}(\alpha_p-\alpha_l)}\right)$$

$$\times\left(\sum_{k=1}^{m-j-3}\sum_{j_1=j+3}^{m-k}\cdots\sum_{j_k=j_{k-1}+1}^{m-1}\frac{\alpha_{j+2\,j_1}\alpha_{j_1 j_2}\cdots\alpha_{j_k m}}{(\alpha_m-\alpha_{j+2})(\alpha_m-\alpha_{j_1})\cdots(\alpha_m-\alpha_{j_k})} + \frac{\alpha_{j+2\,m}}{\alpha_m-\alpha_{j+2}}\right) + \cdots +$$

$$+ \left(\alpha_{ji}\frac{e^{\alpha_j\delta}-e^{\alpha_i\delta}}{\alpha_j-\alpha_i} + \sum_{k=1}^{i-1-j}\sum_{j_1=j+1}^{i-k}\cdots\sum_{j_k=j_{k-1}+1}^{i-1}\sum_{p\in S_j^i(\mathbb{J}_k)}\frac{\alpha_{jj_1}\alpha_{j_1 j_2}\cdots\alpha_{j_k i}\ e^{\alpha_p\delta}}{\prod_{l\in S_j^i(\mathbb{J}_k):l\neq p}(\alpha_p-\alpha_l)}\right)$$

$$\times\left(\sum_{k=1}^{m-i-1}\sum_{j_1=i+1}^{m-k}\cdots\sum_{j_k=j_{k-1}+1}^{m-1}\frac{\alpha_{ij_1}\alpha_{j_1 j_2}\cdots\alpha_{j_k m}}{(\alpha_m-\alpha_i)(\alpha_m-\alpha_{j_1})\cdots(\alpha_m-\alpha_{j_k})} + \frac{\alpha_{im}}{\alpha_m-\alpha_i}\right)$$

$$+$$
$$\vdots$$
$$+$$

$$+ \left(\alpha_{jm-2}\frac{e^{\alpha_j\delta}-e^{\alpha_{m-2}\delta}}{\alpha_j-\alpha_{m-2}} + \sum_{k=1}^{m-j-3}\sum_{j_1=j+1}^{m-2-k}\cdots\sum_{j_k=j_{k-1}+1}^{m-3}\sum_{p\in S_j^{m-2}(\mathbb{J}_k)}\frac{\alpha_{jj_1}\alpha_{j_1 j_2}\cdots\alpha_{j_k m-2}e^{\alpha_p\delta}}{\prod_{l\in S_j^{m-2}(\mathbb{J}_k):l\neq p}(\alpha_p-\alpha_l)}\right)$$

$$\times\left(\frac{\alpha_{m-2\,m-1}\alpha_{m-1\,m}}{(\alpha_m-\alpha_{m-2})(\alpha_m-\alpha_{m-1})} + \frac{\alpha_{m-2\ m}}{\alpha_m-\alpha_{m-2}}\right)$$

---

[5]$e^{B_m\delta} = I + B_m\delta + \frac{(B_m\delta)^2}{2!} + + \frac{(B_m\delta)^3}{3!} + \cdots$; where $I$ is the identity matrix of appropriate order.

[6]The first term is $(e^{B\delta})_{jj}a_j^m$, the second term is $(e^{B\delta})_{j\,j+1}a_{j+1}$ Note that equation (49) of Lemma 1 stops with the first term in RHS when $i = p+1$, and $a_m^m = 1$.

$$+ \left( \alpha_{jm-1} \frac{e^{\alpha_j \delta} - e^{\alpha_{m-1} \delta}}{\alpha_j - \alpha_{m-1}} + \sum_{k=1}^{m-j-2} \sum_{j_1=j+1}^{m-1-k} \cdots \sum_{j_k=j_{k-1}+1}^{m-2} \sum_{p \in S_j^{m-1}(\mathbb{J}_k)} \frac{\alpha_{jj_1} \alpha_{j_1 j_2} \cdots \alpha_{j_k m-1} e^{\alpha_p \delta}}{\prod_{l \in S_j^{m-1}(\mathbb{J}_k): l \neq p} (\alpha_p - \alpha_l)} \right) \frac{\alpha_{m-1 m}}{\alpha_m - \alpha_{m-1}} +$$

$$+ \frac{e^{\alpha_j \delta} - e^{\alpha_m \delta}}{\alpha_j - \alpha_m} \alpha_{jm} + \sum_{k=1}^{m-j-1} \sum_{j_1=j+1}^{m-k} \cdots \sum_{j_k=j_{k-1}+1}^{m-1} \sum_{p \in S_j^m(\mathbb{J}_k)} \frac{\alpha_{jj_1} \alpha_{j_1 j_2} \cdots \alpha_{j_k m} e^{\alpha_p \delta}}{\prod_{l \in S_j^m(\mathbb{J}_k): l \neq p} (\alpha_p - \alpha_l)}.$$

*One can write the above in compact form as*

$$\sum_{i=j}^{m} \left( e^{B_m \delta} \right)_{ji} a_i^m = e^{\alpha_j \delta} \left( \sum_{k=1}^{m-j-1} \sum_{j_1=j+1}^{m-k} \cdots \sum_{j_k=j_{k-1}+1}^{m-1} \frac{\alpha_{jj_1} \alpha_{j_1 j_2} \cdots \alpha_{j_k m}}{(\alpha_m - \alpha_j)(\alpha_m - \alpha_{j_1}) \cdots (\alpha_m - \alpha_{j_k})} + \frac{\alpha_{jm}}{\alpha_m - \alpha_j} \right)$$

$$+ \alpha_{jj+1} \frac{e^{\alpha_j \delta} - e^{\alpha_{j+1} \delta}}{\alpha_j - \alpha_{j+1}} \left( \sum_{k=1}^{m-j-2} \sum_{j_1=j+2}^{m-k} \cdots \sum_{j_k=j_{k-1}+1}^{m-1} \frac{\alpha_{j+1 j_1} \alpha_{j_1 j_2} \cdots \alpha_{j_k m}}{(\alpha_m - \alpha_{j+1})(\alpha_m - \alpha_{j_1}) \cdots (\alpha_m - \alpha_{j_k})} + \frac{\alpha_{j+1 m}}{\alpha_m - \alpha_{j+1}} \right)$$

$$+ \sum_{i=j+2}^{m-2} \left( \alpha_{ji} \frac{e^{\alpha_j \delta} - e^{\alpha_i \delta}}{\alpha_j - \alpha_i} + \sum_{k=1}^{i-1-j} \sum_{j_1=j+1}^{i-k} \cdots \sum_{j_k=j_{k-1}+1}^{i-1} \sum_{p \in S_j^i(\mathbb{J}_k)} \frac{\alpha_{jj_1} \alpha_{j_1 j_2} \cdots \alpha_{j_k i} e^{\alpha_p \delta}}{\prod_{l \in S_j^i(\mathbb{J}_k): l \neq p} (\alpha_p - \alpha_l)} \right)$$

$$\times \left( \sum_{k=1}^{m-i-1} \sum_{j_1=i+1}^{m-k} \cdots \sum_{j_k=j_{k-1}+1}^{m-1} \frac{\alpha_{ij_1} \alpha_{j_1 j_2} \cdots \alpha_{j_k m}}{(\alpha_m - \alpha_i)(\alpha_m - \alpha_{j_1}) \cdots (\alpha_m - \alpha_{j_k})} + \frac{\alpha_{im}}{\alpha_m - \alpha_i} \right)$$

$$+ \left( \alpha_{jm-1} \frac{e^{\alpha_j \delta} - e^{\alpha_{m-1} \delta}}{\alpha_j - \alpha_{m-1}} + \sum_{k=1}^{m-j-2} \sum_{j_1=j+1}^{m-1-k} \cdots \sum_{j_k=j_{k-1}+1}^{m-2} \sum_{p \in S_j^{m-1}(\mathbb{J}_k)} \frac{\alpha_{jj_1} \alpha_{j_1 j_2} \cdots \alpha_{j_k m-1} e^{\alpha_p \delta}}{\prod_{l \in S_j^{m-1}(\mathbb{J}_k): l \neq p} (\alpha_p - \alpha_l)} \right) \frac{\alpha_{m-1 m}}{\alpha_m - \alpha_{m-1}} +$$

$$+ \frac{e^{\alpha_j \delta} - e^{\alpha_m \delta}}{\alpha_j - \alpha_m} \alpha_{jm} + \sum_{k=1}^{m-j-1} \sum_{j_1=j+1}^{m-k} \cdots \sum_{j_k=j_{k-1}+1}^{m-1} \sum_{p \in S_j^m(\mathbb{J}_k)} \frac{\alpha_{jj_1} \alpha_{j_1 j_2} \cdots \alpha_{j_k m} e^{\alpha_p \delta}}{\prod_{l \in S_j^m(\mathbb{J}_k): l \neq p} (\alpha_p - \alpha_l)}. \quad (29)$$

*We claim that the simplification of the equation (29) yields the terms containing $e^{\alpha_m \delta}$ only and all the others cancel. We first gather the coefficient of $e^{\alpha_m \delta}$ in the equation (29). Observe that the terms containing $e^{\alpha_m \delta}$ appear only in the last line in this equation with $p = m$ and is given as*

$$\sum_{k=1}^{m-j-1} \sum_{j_1=j+1}^{m-k} \cdots \sum_{j_k=j_{k-1}+1}^{m-1} \frac{\alpha_{jj_1} \alpha_{j_1 j_2} \cdots \alpha_{j_k m} e^{\alpha_m \delta}}{(\alpha_m - \alpha_j)(\alpha_m - \alpha_{j_1}) \cdots (\alpha_m - \alpha_{j_k})} + \frac{\alpha_{jm} e^{\alpha_m \delta}}{\alpha_m - \alpha_j} \quad \text{which is equal to } a_j^m e^{\alpha_m \delta}. \quad (30)$$

Thus, we get $a_j^m$ as the net coefficient of $e^{\alpha_m \delta}$. We will now show the net coefficient of $e^{\alpha_l \delta}$, where $l \neq m$ are all zero. Define $\Upsilon_j$ as the net coefficient of $e^{\alpha_j \delta}$ (in equation (29)) and we begin with $j = m - 1$. Basically $\Upsilon_{m-1}$ is obtained with $p = m - 1$ in the last two lines of equation (29) and one additional term in the second last line.

$$\Upsilon_{m-1} = \left( -\frac{\alpha_{jm-1} e^{\alpha_{m-1} \delta}}{\alpha_j - \alpha_{m-1}} + \sum_{k=1}^{m-j-2} \sum_{j_1=j+1}^{m-1-k} \cdots \sum_{j_k=j_{k-1}+1}^{m-2} \frac{\alpha_{jj_1} \alpha_{j_1 j_2} \cdots \alpha_{j_k m-1} e^{\alpha_{m-1} \delta}}{(\alpha_{m-1} - \alpha_j)(\alpha_{m-1} - \alpha_{j_1}) \cdots (\alpha_{m-1} - \alpha_{j_k})} \right) \frac{\alpha_{m-1 m}}{\alpha_m - \alpha_{m-1}}$$

$$+ \sum_{k=1}^{m-j-1} \sum \sum \cdots \sum \left\{ \begin{array}{l} j_k = m-1: \quad j+1 \leq j_1 \leq m-k, \\ \vdots \\ j_{k-3}+1 \leq j_{k-2} \leq m-3 \\ j_{k-1}+1 \leq m-1 \end{array} \right\} \frac{\alpha_{jj_1} \alpha_{j_1 j_2} \cdots \alpha_{j_k m} e^{\alpha_{m-1} \delta}}{(\alpha_{m-1} - \alpha_j)(\alpha_{m-1} - \alpha_{j_1}) \cdots (\alpha_{m-1} - \alpha_{j_k})}$$

$$
\begin{aligned}
=\ & -\frac{\alpha_{jm-1}\alpha_{m-1m}e^{\alpha_{m-1}\delta}}{(\alpha_{m-1}-\alpha_j)(\alpha_{m-1}-\alpha_m)} \\
& -\sum_{k=1}^{m-j-2}\sum_{j_1=j+1}^{m-1-k}\cdots\sum_{j_k=j_{k-1}+1}^{m-2}\frac{\alpha_{jj_1}\alpha_{j_1j_2}\cdots\alpha_{j_km-1}\alpha_{m-1m}e^{\alpha_{m-1}\delta}}{(\alpha_{m-1}-\alpha_j)(\alpha_{m-1}-\alpha_{j_1})\cdots(\alpha_{m-1}-\alpha_{j_k})(\alpha_{m-1}-\alpha_m)} \\
& +\frac{\alpha_{jm-1}\alpha_{m-1m}e^{\alpha_{m-1}\delta}}{(\alpha_{m-1}-\alpha_j)(\alpha_{m-1}-\alpha_m)} \\
& +\sum_{k=2}^{m-j-1}\sum\sum\cdots\sum_{\left\{\begin{array}{c} j+1\le j_1\le m-k,\\ \vdots\\ j_{k-3}+1\le j_{k-2}\le m-3\\ j_{k-1}\le m-2\end{array}\right\}}\frac{\alpha_{jj_1}\alpha_{j_1j_2}\cdots\alpha_{j_{k-1}m-1}\alpha_{m-1m}e^{\alpha_{m-1}\delta}}{(\alpha_{m-1}-\alpha_j)(\alpha_{m-1}-\alpha_{j_1})\cdots(\alpha_{m-1}-\alpha_{j_{k-1}})} \quad (\text{substituting } j_k)
\end{aligned}
$$

$$
\begin{aligned}
=\ & -\frac{\alpha_{jm-1}\alpha_{m-1m}e^{\alpha_{m-1}\delta}}{(\alpha_{m-1}-\alpha_j)(\alpha_{m-1}-\alpha_m)} \\
& -\sum_{k=1}^{m-j-2}\sum_{j_1=j+1}^{m-1-k}\cdots\sum_{j_k=j_{k-1}+1}^{m-2}\frac{\alpha_{jj_1}\alpha_{j_1j_2}\cdots\alpha_{j_km-1}\alpha_{m-1m}e^{\alpha_{m-1}\delta}}{(\alpha_{m-1}-\alpha_j)(\alpha_{m-1}-\alpha_{j_1})\cdots(\alpha_{m-1}-\alpha_{j_k})(\alpha_{m-1}-\alpha_m)} \\
& +\frac{\alpha_{jm-1}\alpha_{m-1m}e^{\alpha_{m-1}\delta}}{(\alpha_{m-1}-\alpha_j)(\alpha_{m-1}-\alpha_m)} \\
& +\sum_{k=1}^{m-j-2}\sum\sum\cdots\sum_{\left\{\begin{array}{c} j+1\le j_1\le m-k,\\ \vdots\\ j_{k-2}+1\le j_{k-1}\le m-3\\ j_k\le m-2\end{array}\right\}}\frac{\alpha_{jj_1}\alpha_{j_1j_2}\cdots\alpha_{j_km-1}\alpha_{m-1m}e^{\alpha_{m-1}\delta}}{(\alpha_{m-1}-\alpha_j)(\alpha_{m-1}-\alpha_{j_1})\cdots(\alpha_{m-1}-\alpha_{j_k})(\alpha_{m-1}-\alpha_m)}=0
\end{aligned}
$$

Hence we see that $\Upsilon_{m-1}=0$.

In similar way, the coefficient of $e^{\alpha_{m-2}\delta}$ is obtained by taking the last three lines of the equation 29 and so on.

Observe that $\Upsilon_j$ (see equation (29)) can be written as $\Upsilon_j=\sum_l \Upsilon_j^l$, where $\Upsilon_j^l$ are the components containing the product of $l$–terms in the numerator and appearing with $e^{\alpha_j\delta}$. Further, note that $\Upsilon_j$ involves the greatest complexity among all in terms of number of terms appearing. In general, we have

$$Comp(\Upsilon_j)>Comp(\Upsilon_{j+1})\cdots>Comp(\Upsilon_{m-1}),\ \text{where } Comp(\Upsilon_i)\ \text{is the complexity of}\ \Upsilon_i.$$

Thus, the computations of $\Upsilon_{j+1},\cdots\Upsilon_{m-2}$ and so on follow quite similar method as applied on a part of $\Upsilon_j$. One can easily verify that $\Upsilon_{j+1},\Upsilon_{j+2},\cdots$ to be zero once we show $\Upsilon_j=0$.[7] Towards showing $\Upsilon_j=0$, we first gather $\Upsilon_j^1$, i.e., all terms containing $e^{\alpha_j\delta}$ appearing with one term in the numerator (coming from the last term of the first line and the first term of the last line of the equation (29)):

$$e^{\alpha_j\delta}\frac{\alpha_{jm}}{\alpha_m-\alpha_j}+\frac{e^{\alpha_j\delta}}{\alpha_j-\alpha_m}\alpha_{jm}=0.$$

Now gather all terms containing $e^{\alpha_j\delta}$ appearing with product of two terms in the numerator. These are obtained by putting $k=1$ in the first and the last line (additionally $p=j$ for the last line)[8], and collecting the

---

[7] Illustration for computation $\Upsilon_{j+1}$ is given after showing $\Upsilon_j=0$ and rest of them follow similar logic.
[8] Note that $\mathbb{J}_1=\{j_1\}$ and $S_j^m(\mathbb{J}_1)=\{j,j_1,m\}$.

*intermediate terms which contain $\alpha_{jj+1}\alpha_{j+1m}, \cdots$ so on as their parts. So,*

$$\begin{aligned}
\Upsilon_j^2 &= \sum_{j_1=j+1}^{m-1} \frac{\alpha_{jj_1}\alpha_{j_1m}e^{\alpha_j\delta}}{(\alpha_m-\alpha_j)(\alpha_m-\alpha_{j_1})} - \alpha_{jj+1}\frac{\alpha_{j+1m}e^{\alpha_j\delta}}{(\alpha_j-\alpha_{j+1})(\alpha_{j+1}-\alpha_m)} - \alpha_{jj+2}\frac{\alpha_{j+2m}e^{\alpha_j\delta}}{(\alpha_j-\alpha_{j+2})(\alpha_{j+2}-\alpha_m)} - \cdots \\
&\quad - \cdots - \frac{\alpha_{jm-1}\alpha_{m-1m}e^{\alpha_j\delta}}{(\alpha_j-\alpha_{m-1})(\alpha_{m-1}-\alpha_m)} + \sum_{j_1=j+1}^{m-1}\frac{\alpha_{jj_1}\alpha_{j_1m}e^{\alpha_j\delta}}{(\alpha_j-\alpha_{j_1})(\alpha_j-\alpha_m)} \\
&= \sum_{j_1=j+1}^{m-1}\frac{\alpha_{jj_1}\alpha_{j_1m}e^{\alpha_j\delta}}{(\alpha_j-\alpha_m)(\alpha_{j_1}-\alpha_m)} - \sum_{j_1=j+1}^{m-1}\frac{\alpha_{jj_1}\alpha_{j_1m}e^{\alpha_j\delta}}{(\alpha_j-\alpha_{j_1})(\alpha_{j_1}-\alpha_m)} + \sum_{j_1=j+1}^{m-1}\frac{\alpha_{jj_1}\alpha_{j_1m}e^{\alpha_j\delta}}{(\alpha_j-\alpha_m)(\alpha_j-\alpha_{j_1})} \\
&= \sum_{j_1=j+1}^{m-1}\frac{\alpha_{jj_1}\alpha_{j_1m}e^{\alpha_j\delta}}{(\alpha_{j_1}-\alpha_m)}\left(\frac{1}{(\alpha_j-\alpha_m)} - \frac{1}{\alpha_j-\alpha_{j_1}}\right) + \sum_{j_1=j+1}^{m-1}\frac{\alpha_{jj_1}\alpha_{j_1m}e^{\alpha_j\delta}}{(\alpha_j-\alpha_{j_1})(\alpha_j-\alpha_m)} \\
&= \sum_{j_1=j+1}^{m-1}\frac{\alpha_{jj_1}\alpha_{j_1m}e^{\alpha_j\delta}}{(\alpha_{j_1}-\alpha_m)}\frac{\alpha_m-\alpha_{j_1}}{(\alpha_j-\alpha_m)(\alpha_j-\alpha_{j_1})} + \sum_{j_1=j+1}^{m-1}\frac{\alpha_{jj_1}\alpha_{j_1m}e^{\alpha_j\delta}}{(\alpha_j-\alpha_m)(\alpha_j-\alpha_{j_1})} = 0.
\end{aligned}$$

*Thus, we get zero as the coefficient of $e^{\alpha_j\delta}$ when it appears with the product of the two terms, i.e., $\Upsilon_j^2 = 0$. We will now consider the overall coefficient of $e^{\alpha_j\delta}$ constructed using all components have $r \geq 3$ number of terms in the numerator $\Upsilon_j^r$ and show that the overall coefficient $\Upsilon_j^r$ is again zero. Once this is done, the proof is complete for $j$, i.e., the sum of all the terms containing term $e^{\alpha_j\delta}$, $\Upsilon_j = 0$.*

*Towards this we will have terms like the following, each of which can be split as below into 4 different sub-terms (one for each i and note only terms containing $e^{\alpha_j\delta}$ (with $p = j$ in the second one) are considered); and from among these we need to extract the ones that can be a part of the computation of $\Upsilon_j^r$:*

$$\begin{aligned}
&\left(\frac{\alpha_{ji}e^{\alpha_j\delta}}{\alpha_j-\alpha_i} + \sum_{k=1}^{i-1-j}\sum_{j_1=j+1}^{i-k}\cdots\sum_{j_k=j_{k-1}+1}^{i-1}\frac{\alpha_{jj_1}\alpha_{j_1j_2}\cdots\alpha_{j_ki}e^{\alpha_j\delta}}{\prod_{l\in S_j^i(\mathbb{J}_k):l\neq j}(\alpha_j-\alpha_l)}\right) \\
&\qquad\times\left(\sum_{k=1}^{m-i-1}\sum_{j_1=i+1}^{m-k}\cdots\sum_{j_k=j_{k-1}+1}^{m-1}\frac{\alpha_{ij_1}\alpha_{j_1j_2}\cdots\alpha_{j_km}}{\prod_{l\in S_i^m(\mathbb{J}_k):l\neq m}(\alpha_m-\alpha_l)} + \frac{\alpha_{im}}{\alpha_m-\alpha_i}\right) \\
&= \sum_{k=1}^{m-i-1}\sum_{j_1=i+1}^{m-k}\cdots\sum_{j_k=j_{k-1}+1}^{m-1}\frac{\alpha_{ji}e^{\alpha_j\delta}}{\alpha_j-\alpha_i}\frac{\alpha_{ij_1}\alpha_{j_1j_2}\cdots\alpha_{j_km}}{\prod_{l\in S_i^m(\mathbb{J}_k):l\neq m}(\alpha_m-\alpha_l)} \\
&\quad + \sum_{k=1}^{i-1-j}\sum_{j_1=j+1}^{i-k}\cdots\sum_{j_k=j_{k-1}+1}^{i-1}\frac{\alpha_{jj_1}\alpha_{j_1j_2}\cdots\alpha_{j_ki}e^{\alpha_j\delta}}{\prod_{l\in S_j^i(\mathbb{J}_k):l\neq j}(\alpha_j-\alpha_l)}\frac{\alpha_{im}}{\alpha_m-\alpha_i} \qquad (31) \\
&\quad + \sum_{k=1}^{i-1-j}\sum_{j_1=j+1}^{i-k}\cdots\sum_{j_k=j_{k-1}+1}^{i-1}\frac{\alpha_{jj_1}\alpha_{j_1j_2}\cdots\alpha_{j_ki}e^{\alpha_j\delta}}{\prod_{l\in S_j^i(\mathbb{J}_k):l\neq j}(\alpha_j-\alpha_l)} \\
&\qquad\times\sum_{k'=1}^{m-i-1}\sum_{j_1=i+1}^{m-k'}\cdots\sum_{j'_k=j_{k'-1}+1}^{m-1}\frac{\alpha_{ij_1}\alpha_{j_1j_2}\cdots\alpha_{j_{k'}m}}{\prod_{l\in S_i^m(\mathbb{J}_{k'}):l\neq m}(\alpha_m-\alpha_l)} + \frac{\alpha_{ji}\alpha_{im}}{(\alpha_j-\alpha_i)(\alpha_m-\alpha_i)} \qquad (32)
\end{aligned}$$

*Using the above split after gathering required terms with product of $r$-terms in the numerator (from (29 and*

the Lemma) we have[9]:

$$
\begin{aligned}
\Upsilon_j^r &= \underbrace{\sum_{j_1=j+1}^{m-r+1} \sum_{j_2=j_1+1}^{m-r+2} \cdots \sum_{j_{r-1}=j_{r-2}+1}^{m-1} \frac{\alpha_{jj_1}\alpha_{j_1j_2}\cdots\alpha_{j_{r-1}m}e^{\alpha_j\delta}}{(\alpha_m-\alpha_j)(\alpha_m-\alpha_{j_1})\cdots(\alpha_m-\alpha_{j_{r-1}})}}_{\text{from }(e^{B_m\delta})_{j,j}a_j^m \text{ component}} \\
&+ \underbrace{\sum_{i=j+1}^{m-r+1} \sum_{j_1=i+1}^{m-r+2} \cdots \sum_{j_{r-2}=j_{r-3}+1}^{m-1} \frac{\alpha_{ji}\alpha_{ij_1}\alpha_{j_1j_2}\cdots\alpha_{j_{r-2}m}e^{\alpha_j\delta}}{(\alpha_j-\alpha_i)(\alpha_m-\alpha_i)(\alpha_m-\alpha_{j_1})(\alpha_m-\alpha_{j_2})\cdots(\alpha_m-\alpha_{j_{r-2}})}}_{\substack{\text{from}(e^{B_m\delta})_{j,i}a_i^m \text{ components, with } i \text{ s.t. it becomes product of } r \text{ terms } (k=r-2) \\ \text{and using first part of the split equation (31)}}}
\end{aligned}
$$

$$
+ \underbrace{\sum_{j_1=j+1}^{m-r+1} \sum_{j_2=j_1+1}^{m-r+2} \cdots \sum_{j_{r-2}=j_{r-3}+1}^{m-2} \sum_{i=j_{r-2}+1}^{m-1} \frac{\alpha_{jj_1}\alpha_{j_1j_2}\cdots\alpha_{j_{r-2}i}\alpha_{im}e^{\alpha_j\delta}}{(\alpha_j-\alpha_{j_1})\cdots(\alpha_j-\alpha_{j_{r-2}})(\alpha_j-\alpha_i)(\alpha_m-\alpha_i)}}_{\substack{\text{from the second part of the split equation (31) after rearranging} \\ \text{as in footnote 9 and with } k=(r-2) \text{ (s.t. we have product of } r \text{ terms),} \\ \text{i.e., first summation till } j_1 \leq m-r \text{ and one more }(j_1=m-r+1) \text{ term from } (e^{B_m\delta})_{jm-1}a_{m-1}^m}}
$$

$$
+ \sum_{k=1}^{r-3} \sum_{j_1=j+1}^{m-r+1} \cdots \sum_{j_k=j_{k-1}+1}^{m-r+k} \sum_{i=j_k+1}^{m-r+1+k} \frac{\alpha_{jj_1}\alpha_{j_1j_2}\cdots\alpha_{j_ki}}{(\alpha_j-\alpha_{j_1})\cdots(\alpha_j-\alpha_{j_k})(\alpha_j-\alpha_i)}
$$
$$
\times \underbrace{\sum_{j_{k+1}=i+1}^{m-r+k+2} \sum_{j_{k+2}=j_{k+1}+1}^{m-r+k+3} \cdots \sum_{j_{r-2}=j_{r-3}+1}^{m-1} \frac{\alpha_{ij_{k+1}}\alpha_{j_{k+1}j_{k+2}}\cdots\alpha_{j_{r-2}m}}{(\alpha_m-\alpha_i)\left(\alpha_m-\alpha_{j_{k+1}}\right)\cdots\left(\alpha_m-\alpha_{j_{r-2}}\right)}}_{\substack{\text{from the third part of the split equation given by (32)} \\ \text{after rearranging the order of summations as in footnote 9} \\ \text{and using } k'=r-k-2}}
$$

$$
+ \underbrace{\sum_{j_1=j+1}^{m-r+1} \sum_{j_2=j_1+1}^{m-r+2} \cdots \sum_{j_{r-1}=j_{r-2}+1}^{m-1} \frac{\alpha_{jj_1}\alpha_{j_1j_2}\cdots\alpha_{j_{r-1}m}e^{\alpha_j\delta}}{(\alpha_j-\alpha_{j_1})\cdots(\alpha_j-\alpha_{j_{r-1}})(\alpha_j-\alpha_m)}}_{\text{from }(e^{B_m\delta})_{jm}a_m^m} .
$$

---

[9] The order is interchanged from the following (with the convention that $\sum_{j=a}^{b} \cdot\cdot = 0$ when $b < a$)

$$
\begin{aligned}
\sum_{i=j+2}^{m-2} \sum_{k=1}^{i-1-j} \sum_{j_1=j+1}^{i-k} \cdots \sum_{j_k=j_{k-1}+1}^{i-1} &= \sum_{k=1}^{m-3-j} \sum_{i=k+1+j}^{m-2} \sum_{j_1=j+1}^{i-k} \cdots \sum_{j_k=j_{k-1}+1}^{i-1} \\
&= \sum_{k=1}^{m-3-j} \sum_{j_1=j+1}^{m-2-k} \sum_{i=k+j_1}^{m-2} \cdots \sum_{j_k=j_{k-1}+1}^{i-1} \\
&\vdots \\
&= \sum_{k=1}^{m-3-j} \sum_{j_1=j+1}^{m-2-k} \cdots \sum_{j_{k-1}=j_{k-2}+1}^{m-4} \sum_{i=j_{k-1}+2}^{m-2} \sum_{j_k=j_{k-1}+1}^{i-1} \\
&= \sum_{k=1}^{m-3-j} \sum_{j_1=j+1}^{m-2-k} \cdots \sum_{j_{k-1}=j_{k-2}+1}^{m-4} \sum_{j_k=j_{k-1}+1}^{m-3} \sum_{i=j_k+1}^{m-2}
\end{aligned}
$$

*Thus, after using matching indices for all the terms:*

$$\Upsilon_j^r = \sum_{j_1=j+1}^{m-r+1}\sum_{j_2=j_1+1}^{m-r+2}\cdots\sum_{j_{r-1}=j_{r-2}+1}^{m-1}\frac{\alpha_{jj_1}\alpha_{j_1j_2}\cdots\alpha_{j_{r-1}m}e^{\alpha_j\delta}}{(\alpha_m-\alpha_j)(\alpha_m-\alpha_{j_1})\cdots(\alpha_m-\alpha_{j_{r-1}})}$$

$$+\sum_{j_1=j+1}^{m-r+1}\sum_{j_2=j_1+1}^{m-r+2}\cdots\sum_{j_{r-1}=j_{r-2}+1}^{m-1}\frac{\alpha_{jj_1}\alpha_{j_1j_2}\cdots\alpha_{j_{r-1}m}e^{\alpha_j\delta}}{(\alpha_j-\alpha_{j_1})(\alpha_m-\alpha_{j_1})\cdots(\alpha_m-\alpha_{j_{r-1}})}$$

$$+\sum_{j_1=j+1}^{m-r+1}\sum_{j_2=j_1+1}^{m-r+2}\cdots\sum_{j_{r-1}=j_{r-2}+1}^{m-1}\frac{\alpha_{jj_1}\alpha_{j_1j_2}\cdots\alpha_{j_{r-1}m}e^{\alpha_j\delta}}{(\alpha_j-\alpha_{j_1})\cdots(\alpha_j-\alpha_{j_{r-2}})(\alpha_j-\alpha_{j_{r-1}})(\alpha_m-\alpha_{j_{r-1}})}$$

$$+\sum_{k=1}^{r-3}\sum_{j_1=j+1}^{m-r+1}\cdots\sum_{j_k=j_{k-1}+1}^{m-r+k}\sum_{j_{k+1}=j_k+1}^{m-r+1+k}\sum_{j_{k+2}=j_{k+1}+1}^{m-r+k+2}\cdots\sum_{j_{r-1}=j_{r-2}+1}^{m-1}\frac{\alpha_{jj_1}\alpha_{j_1j_2}\cdots\alpha_{j_kj_{k+1}}\alpha_{j_{k+1}j_{k+2}}\cdots\alpha_{j_{r-1}m}}{(\alpha_j-\alpha_{j_1})\cdots(\alpha_j-\alpha_{j_k})(\alpha_j-\alpha_{j_{k+1}})(\alpha_m-\alpha_{j_{k+1}})\cdots(\alpha_m-\alpha_{j_{r-1}})}$$

$$+\sum_{j_1=j+1}^{m-r+1}\sum_{j_2=j_1+1}^{m-r+2}\cdots\sum_{j_{r-1}=j_{r-2}+1}^{m-1}\frac{\alpha_{jj_1}\alpha_{j_1j_2}\cdots\alpha_{j_{r-1}m}e^{\alpha_j\delta}}{(\alpha_j-\alpha_{j_1})\cdots(\alpha_j-\alpha_{j_{r-1}})(\alpha_j-\alpha_m)}.$$

*One can rewrite the above as*

$$\Upsilon_j^r = \sum_{j_1=j+1}^{m-r+1}\sum_{j_2=j_1+1}^{m-r+2}\cdots\sum_{j_{r-1}=j_{r-2}+1}^{m-1}\left(\frac{1}{(\alpha_m-\alpha_j)(\alpha_m-\alpha_{j_1})\cdots(\alpha_m-\alpha_{j_{r-1}})}\right.$$

$$+\frac{1}{(\alpha_j-\alpha_{j_1})(\alpha_m-\alpha_{j_1})(\alpha_m-\alpha_{j_2})\cdots(\alpha_m-\alpha_{j_{r-2}})}$$

$$+\frac{1}{(\alpha_j-\alpha_{j_1})(\alpha_j-\alpha_{j_2})(\alpha_m-\alpha_{j_2})\cdots(\alpha_m-\alpha_{j_{r-2}})}$$

$$+\frac{1}{(\alpha_j-\alpha_{j_1})(\alpha_j-\alpha_{j_2})(\alpha_j-\alpha_{j_3})(\alpha_m-\alpha_{j_3})\cdots(\alpha_m-\alpha_{j_{r-2}})}$$

$$+\cdots+\frac{1}{(\alpha_j-\alpha_{j_1})\cdots(\alpha_j-\alpha_{j_{r-2}})(\alpha_m-\alpha_{j_{r-2}})(\alpha_m-\alpha_{j_{r-1}})}$$

$$\left.+\frac{1}{(\alpha_j-\alpha_{j_1})\cdots(\alpha_j-\alpha_{j_{r-1}})(\alpha_j-\alpha_m)}\right). \tag{33}$$

*and observe the following for any $(j_1, j_2, \cdots j_{r-1})$:*

$$\underbrace{\frac{1}{(\alpha_m-\alpha_j)(\alpha_m-\alpha_{j_1})\cdots(\alpha_m-\alpha_{j_{r-1}})}+\frac{1}{(\alpha_j-\alpha_{j_1})(\alpha_m-\alpha_{j_1})(\alpha_m-\alpha_{j_2})\cdots(\alpha_m-\alpha_{j_{r-2}})}}$$

$$+\frac{1}{(\alpha_j-\alpha_{j_1})(\alpha_j-\alpha_{j_2})(\alpha_m-\alpha_{j_2})\cdots(\alpha_m-\alpha_{j_{r-2}})}$$

$$+\frac{1}{(\alpha_j-\alpha_{j_1})(\alpha_j-\alpha_{j_2})(\alpha_j-\alpha_{j_3})(\alpha_m-\alpha_{j_3})\cdots(\alpha_m-\alpha_{j_{r-2}})}+\cdots+$$

$$+\frac{1}{(\alpha_j-\alpha_{j_1})\cdots(\alpha_j-\alpha_{j_{r-2}})(\alpha_m-\alpha_{j_{r-2}})(\alpha_m-\alpha_{j_{r-1}})}$$

$$+\frac{1}{(\alpha_j-\alpha_{j_1})\cdots(\alpha_j-\alpha_{j_{r-1}})(\alpha_j-\alpha_m)}$$

$$= \frac{1}{(\alpha_m - \alpha_j)(\alpha_j - \alpha_{j_1})(\alpha_m - \alpha_{j_2})\cdots(\alpha_m - \alpha_{j_{r-2}})} + \frac{1}{(\alpha_j - \alpha_{j_1})(\alpha_j - \alpha_{j_2})(\alpha_m - \alpha_{j_2})\cdots(\alpha_m - \alpha_{j_{r-2}})}$$
$$+ \frac{1}{(\alpha_j - \alpha_{j_1})(\alpha_j - \alpha_{j_2})(\alpha_j - \alpha_{j_3})(\alpha_m - \alpha_{j_3})\cdots(\alpha_m - \alpha_{j_{r-2}})} + \cdots +$$
$$+ \frac{1}{(\alpha_j - \alpha_{j_1})\cdots(\alpha_j - \alpha_{j_{r-2}})(\alpha_m - \alpha_{j_{r-2}})(\alpha_m - \alpha_{j_{r-1}})}$$
$$+ \frac{1}{(\alpha_j - \alpha_{j_1})\cdots(\alpha_j - \alpha_{j_{r-1}})(\alpha_j - \alpha_m)} = 0 \; \forall r \geq 3. \tag{34}$$

Using (34), we get $\Upsilon_j^r = 0$ for all $r \geq 3$. In all, we see that the net coefficient $e^{\alpha_j \delta}$ is zero.

*Illustration for the computation of* $\Upsilon_{j+1}$: It is obtained by substituting $j_1 = j+1$ and $p = j+1$ in the terms coming from $\left(e^{B\delta}\right)_{im}$ and then gather all the terms containing $e^{\alpha_{j+1}\delta}$ using similar method as employed in obtaining $\Upsilon_j$. And so on for the others. Thus the claim is true and hence first result of part (i) follows.

*Proof of (i) (b)*  In this part, we compute the expected value of $X_m(t)$ using equation (25) with $a_i^j$ as defined previously (e.g., see equation (26)), that is:

$$a_i^j = \sum_{k=1}^{j-i-1} \sum_{j_1=i+1}^{j-k} \cdots \sum_{j_k = j_{k-1}+1}^{j-1} \frac{\alpha_{ij_1}\alpha_{j_1 j_2}\cdots \alpha_{j_k j}}{(\alpha_j - \alpha_i)(\alpha_j - \alpha_{j_2})\cdots(\alpha_j - \alpha_{j_k})} + \frac{\alpha_{ij}}{\alpha_j - \alpha_i}; \text{ and } a_i^i = 1 \; \forall i.$$

We have $E[M_m(t)] = E[M_m(0)] = a_1^m$ (using martingale property), and thus

$$E[M_m(t)] = E\left[X_m(t)e^{-\alpha_m t}\right] + \sum_{i=1}^{m-1} a_i^m E\left[X_i(t)e^{-\alpha_m t}\right] = a_1^m e^{\alpha_m t} - \sum_{i=1}^{m-1} a_i^m E[X_i(t)]. \tag{35}$$

**We claim that**

$$E[X_N(t)] = a_1^N e^{\alpha_N t} + \sum_{j=1}^{N-1} P_j^N e^{\alpha_j t}; \qquad \text{where} \tag{36}$$

$$P_1^N = -a_1^N - \sum_{k=1}^{N-2}(-1)^k \sum_{j_1=2}^{N-k} \sum_{j_2=j_1+1}^{N-k+1} \cdots \sum_{j_k = j_{k-1}+1}^{N-1} a_1^{j_1} a_{j_1}^{j_2} \cdots a_{j_k}^N; \quad P_{N-1}^N = -a_1^{N-1} a_{N-1}^N; \text{ and} \tag{37}$$

$$P_j^N = \sum_{k=1}^{N-j-1}(-1)^k \sum_{j_1=j+1}^{N-k} \cdots \sum_{j_k = j_{k-1}+1}^{N-1} a_1^j a_j^{j_1} a_{j_1}^{j_2} \cdots a_{j_k}^N - a_1^j a_j^N, \text{ for } j = 2,3\cdots, N-2. \tag{38}$$

*We use the mathematical induction to prove the claim as follows. For the base case $N = 2$, we have*

$$E[X_2(t)] = a_1^2 e^{\alpha_2 t} - a_1^2 E[X_1(t)] = a_1^2 e^{\alpha_2 t} - a_1^2 e^{\alpha_1 t} = a_1^2 e^{\alpha_2 t} + P_1^2 e^{\alpha_1 t}.$$

The above is because $X_1(t)$ is a standard (non-decomposable) branching process and $E[X_1(t)] = e^{\alpha_1 t}$. Thus, the claim holds for $N = 2$.

We assume that the claim holds true for $N = 3, 4, \cdots, m$, and we prove the same for $N = m+1$ below. Using equation (35), we have

$$E[X_{m+1}(t)] = a_1^{m+1} e^{\alpha_{m+1} t} - a_1^{m+1} E[X_1(t)] - a_2^{m+1} E[X_2(t)] - \cdots - a_m^{m+1} E[X_m(t)].$$

*We now substitute value of* $E[X_j(t)]$ *for* $j = 1, 2, \cdots, m$ *as per the claim (see 36)*

$$E[X_{m+1}(t)] = a_1^{m+1} e^{\alpha_{m+1} t} - a_1^{m+1} e^{\alpha_1 t} - \sum_{p=2}^{m} a_1^p a_p^{m+1} e^{\alpha_p t} - \sum_{p=2}^{m} \sum_{j=1}^{p-1} P_j^p a_p^{m+1} e^{\alpha_j t}; \text{ with } P_i^i = 1 \; \forall i.$$

*Substituting $P_j^p$ using equations (37) and (38) in the above equation, we get*

$$\begin{aligned}
E[X_{m+1}(t)] &= a_1^{m+1} e^{\alpha_{m+1} t} - a_1^{m+1} e^{\alpha_1 t} - \sum_{p=2}^{m} a_1^p a_p^{m+1} e^{\alpha_p t} \\
&+ \sum_{p=2}^{m} a_p^{m+1} \left( a_1^p + \sum_{k=1}^{p-2} (-1)^k \sum_{j_1=2}^{p-k} \sum_{j_2=j_1+1}^{p-k+1} \cdots \sum_{j_k=j_{k-1}+1}^{p-1} a_1^{j_1} a_{j_1}^{j_2} \cdots a_{j_k}^{p} \right) e^{\alpha_1 t} \\
&+ \sum_{p=3}^{m} a_p^{m+1} \left( \sum_{j=2}^{p-1} a_1^j a_j^p e^{\alpha_j t} - \sum_{j=2}^{p-2} \sum_{k=1}^{p-j-1} (-1)^k \sum_{j_1=j+1}^{p-k} \sum_{j_2=j_1+1}^{p-k+1} \cdots \sum_{j_k=j_{k-1}+1}^{p-1} a_1^{j} a_{j}^{j_1} a_{j_1}^{j_2} \cdots a_{j_k}^{p} e^{\alpha_j t} \right) \\
&= a_1^{m+1} \left( e^{\alpha_{m+1} t} - e^{\alpha_1 t} \right) - \sum_{p=2}^{m} a_1^p a_p^{m+1} e^{\alpha_p t} + \sum_{p=2}^{m} a_1^p a_p^{m+1} e^{\alpha_1 t} \\
&+ \sum_{p=2}^{m} \sum_{k=1}^{p-2} (-1)^k \sum_{j_1=2}^{p-k} \sum_{j_2=j_1+1}^{p-k+1} \cdots \sum_{j_k=j_{k-1}+1}^{p-1} a_1^{j_1} a_{j_1}^{j_2} \cdots a_{j_k}^{p} a_p^{m+1} e^{\alpha_1 t} \\
&+ \sum_{p=3}^{m} \sum_{j=2}^{p-1} a_1^j a_j^p a_p^{m+1} e^{\alpha_j t} - \sum_{p=3}^{m} \sum_{j=2}^{p-2} \sum_{k=1}^{p-j-1} (-1)^k \sum_{j_1=j+1}^{p-k} \sum_{j_2=j_1+1}^{p-k+1} \cdots \sum_{j_k=j_{k-1}+1}^{p-1} a_1^{j} a_j^{j_1} a_{j_1}^{j_2} \cdots a_{j_k}^{p} a_p^{m+1} e^{\alpha_j t}. \quad (39)
\end{aligned}$$

*Changing the order of summations:*

$$\begin{aligned}
&\sum_{p=3}^{m} \sum_{k=1}^{p-2} (-1)^k \sum_{j_1=2}^{p-k} \sum_{j_2=j_1+1}^{p-k+1} \cdots \sum_{j_k=j_{k-1}+1}^{p-1} a_1^{j_1} a_{j_1}^{j_2} \cdots a_{j_k}^{p} a_p^{m+1} \\
&= \sum_{k=1}^{m-2} (-1)^k \sum_{p=k+2}^{m} \sum_{j_1=2}^{p-k} \sum_{j_2=j_1+1}^{p-k+1} \cdots \sum_{j_k=j_{k-1}+1}^{p-1} a_1^{j_1} a_{j_1}^{j_2} \cdots a_{j_k}^{p} a_p^{m+1} \\
&= \sum_{k=1}^{m-2} (-1)^k \sum_{j_1=2}^{m-k} \sum_{p=j_1+k}^{m} \sum_{j_2=j_1+1}^{p-k+1} \cdots \sum_{j_k=j_{k-1}+1}^{p-1} a_1^{j_1} a_{j_1}^{j_2} \cdots a_{j_k}^{p} a_p^{m+1} \\
&= \sum_{k=1}^{m-2} (-1)^k \sum_{j_1=2}^{m-k} \sum_{j_2=j_1+1}^{m-k+1} \sum_{p=j_2+k-1}^{m} \cdots \sum_{j_k=j_{k-1}+1}^{p-1} a_1^{j_1} a_{j_1}^{j_2} \cdots a_{j_k}^{p} a_p^{m+1} \\
&= \ldots = \sum_{k=1}^{m-2} (-1)^k \sum_{j_1=2}^{m-k} \sum_{j_2=j_1+1}^{m-k+1} \cdots \sum_{j_k=j_{k-1}+1}^{m-1} \sum_{p=j_k+1}^{m} a_1^{j_1} a_{j_1}^{j_2} \cdots a_{j_k}^{p} a_p^{m+1}.
\end{aligned}$$

*To summarize (changing order), we have*

$$\begin{aligned}
&\sum_{p=3}^{m} \sum_{k=1}^{p-2} (-1)^k \sum_{j_1=2}^{p-k} \sum_{j_2=j_1+1}^{p-k+1} \cdots \sum_{j_k=j_{k-1}+1}^{p-1} a_1^{j_1} a_{j_1}^{j_2} \cdots a_{j_k}^{p} a_p^{m+1} \\
&= \sum_{k=1}^{m-2} (-1)^k \sum_{j_1=2}^{m-k} \sum_{j_2=j_1+1}^{m-k+1} \cdots \sum_{j_k=j_{k-1}+1}^{m-1} \sum_{p=j_k+1}^{m} a_1^{j_1} a_{j_1}^{j_2} \cdots a_{j_k}^{p} a_p^{m+1}. \quad (40)
\end{aligned}$$

*Observe that*

$$\begin{aligned}
&\sum_{p=2}^{m} a_1^p a_p^{m+1} + \sum_{k=1}^{m-2} (-1)^k \sum_{j_1=2}^{m-k} \sum_{j_2=j_1+1}^{m-k+1} \cdots \sum_{j_k=j_{k-1}+1}^{m-1} \sum_{p=j_k+1}^{m} a_1^{j_1} a_{j_1}^{j_2} \cdots a_{j_k}^{p} a_p^{m+1} \\
&= -\sum_{k=1}^{m-1} (-1)^k \sum_{j_1=2}^{m+1-k} \sum_{j_2=j_1+1}^{m-k+2} \cdots \sum_{j_k=j_{k-1}+1}^{m} a_1^{j_1} a_{j_1}^{j_2} \cdots a_{j_k}^{m+1}; \\
\Longrightarrow &\sum_{p=2}^{m} a_1^p a_p^{m+1} + \sum_{p=3}^{m} \sum_{k=1}^{p-2} (-1)^k \sum_{j_1=2}^{p-k} \sum_{j_2=j_1+1}^{p-k+1} \cdots \sum_{j_k=j_{k-1}+1}^{p-1} a_1^{j_1} a_{j_1}^{j_2} \cdots a_{j_k}^{p} a_p^{m+1} \\
&= -\sum_{k=1}^{m-1} (-1)^k \sum_{j_1=2}^{m+1-k} \sum_{j_2=j_1+1}^{m-k+2} \cdots \sum_{j_k=j_{k-1}+1}^{m} a_1^{j_1} a_{j_1}^{j_2} \cdots a_{j_k}^{m+1}. \quad (41)
\end{aligned}$$

*Similarly, it is easy to verify that*

$$-\sum_{p=2}^{m} a_1^p a_p^{m+1} e^{\alpha_p t} + \sum_{p=3}^{m} \sum_{j=2}^{p-1} a_1^j a_j^p a_p^{m+1} e^{\alpha_j t}$$

$$-\sum_{p=3}^{m} \sum_{j=2}^{p-2} \sum_{k=1}^{p-j-1} (-1)^k \sum_{j_1=j+1}^{p-k} \sum_{j_2=j_1+1}^{p-k+1} \cdots \sum_{j_k=j_{k-1}+1}^{p-1} a_1^j a_j^{j_1} a_{j_1}^{j_2} \cdots a_{j_k}^p a_p^{m+1} e^{\alpha_j t}$$

$$= -\sum_{j=2}^{m} a_1^j a_j^{m+1} e^{\alpha_j t} + \sum_{j=2}^{m-1} \sum_{k=1}^{m-j} (-1)^k \sum_{j_1=j+1}^{m-k+1} \cdots \sum_{j_k=j_{k-1}+1}^{m} a_1^j a_j^{j_1} a_{j_1}^{j_2} \cdots a_{j_k}^{m+1} e^{\alpha_j t} \qquad (42)$$

*Substituting the expressions of equations (41) and (42) in 39, we obtain*

$$\begin{aligned}
E[X_{m+1}(t)] &= a_1^{m+1}\left(e^{\alpha_{m+1} t} - e^{\alpha_1 t}\right) - \sum_{k=1}^{m-1} (-1)^k \sum_{j_1=2}^{m+1-k} \sum_{j_2=j_1+1}^{m-k+2} \cdots \sum_{j_k=j_{k-1}+1}^{m} a_1^{j_1} a_{j_1}^{j_2} \cdots a_{j_k}^{m+1} \\
&\quad - \sum_{j=2}^{m} a_1^j a_j^{m+1} e^{\alpha_j t} + \sum_{j=2}^{m-1} \sum_{k=1}^{m-j} (-1)^k \sum_{j_1=j+1}^{m-k+1} \cdots \sum_{j_k=j_{k-1}+1}^{m} a_1^j a_j^{j_1} a_{j_1}^{j_2} \cdots a_{j_k}^{m+1} e^{\alpha_j t} \\
&= a_1^{m+1} e^{\alpha_{m+1} t} + \sum_{j=1}^{m} P_j^{m+1} e^{\alpha_j t}; \qquad \text{and this completes the induction based proof.}
\end{aligned}$$

*Hence we have $E[X_m(t)] = a_1^m e^{\alpha_m t} + \sum_{j=1}^{m-1} P_j^m e^{\alpha_j t}$.*

**Proof of (ii)** *When $\alpha_1 > \alpha_2 > \cdots > \alpha_m$, the martingale $M_m(t)$ is positive for all $t \geq 0$. Using the fact that $M_m(t)$ is positive, it converges almost surely to non-negative integrable random variable $W_m$ (e.g., see [12]). And hence*

$$X_m(t) e^{-\alpha_m t} + \sum_{k=1}^{m-1} \sum_{j_1=1}^{m-k} \sum_{j_2=j_1+1}^{m-k+1} \sum_{j_3=j_2+1}^{m-k+2} \cdots \sum_{j_k=j_{k-1}+1}^{m-1} \frac{\alpha_{j_1 j_2} \alpha_{j_2 j_3} \cdots \alpha_{j_k m} X_{j_1}(t) e^{-\alpha_m t}}{(\alpha_m - \alpha_{j_1})(\alpha_m - \alpha_{j_2}) \cdots (\alpha_m - \alpha_{j_k})} \xrightarrow{a.s.} W_m$$

*Equivalently,* $\left| X_m(t) e^{-\alpha_m t} + \sum_{i=1}^{m-1} a_i^m X_i(t) e^{-\alpha_m t} - W_m \right| \xrightarrow{a.s.} 0.$ \hfill (43)

*We show the following with $P_j^m$ (as defined previously).*

$$\left| X_m(t) e^{-\alpha_m t} - W_1 P_1^m e^{\alpha_1 t} e^{-\alpha_m t} + \sum_{j=2}^{m-1} W_j a_j^m e^{(\alpha_j - \alpha_m) t} \right.$$

$$\left. - \sum_{j=2}^{m-2} \sum_{k=1}^{m-j-1} (-1)^k \sum_{j_1=j+1}^{m-k} \cdots \sum_{j_k=j_{k-1}+1}^{m-1} W_j a_j^{j_1} a_{j_1}^{j_2} \cdots a_{j_k}^m e^{(\alpha_j - \alpha_m) t} - W_m \right| \xrightarrow{a.s.} 0.$$

*In a compact form, we claim that*

$$\left| X_N(t) e^{-\alpha_N t} - \sum_{i=1}^{N-1} \frac{P_i^N}{a_1^i} W_i e^{(\alpha_i - \alpha_N) t} - W_N \right| \xrightarrow{a.s.} 0. \qquad (44)$$

*As before, we again use mathematical induction to prove this result. Base case, i.e., for $N = 2$, we have (using equation (43))*

$$\left| X_2(t) e^{-\alpha_2 t} + a_1^2 X_1(t) e^{-\alpha_2 t} - W_2 \right| \xrightarrow{a.s.} 0, \text{ and } \left| W_1 - X_1(t) e^{-\alpha_1 t} \right| \xrightarrow{a.s.} 0$$

*where $W_1, W_2$ are non-negative random variables. The base case immediately follows from the fact that $P_1^2 = -a_1^2$ and $a_1^1 = 1$, hence $\left| X_2(t) e^{-\alpha_2 t} - \frac{P_1^2}{a_1^1} X_1(t) e^{-\alpha_2 t} - W_2 \right| \xrightarrow{a.s.} 0.$*

*Now assume that the result (44 holds for all $N = 3, \cdots, m-1$. We now show it for $N = m$ as follow:*

$$X_m(t) e^{-\alpha_m t} + \sum_{i=1}^{m-1} a_i^m X_i(t) e^{-\alpha_m t} - W_m = X_m(t) e^{-\alpha_m t} - W_m + a_1^m e^{(\alpha_1 - \alpha_m) t} \left( X_1(t) e^{-\alpha_1 t} - W_1 + W_1 \right) +$$

$$+ \sum_{k=2}^{m-1} a_k^m e^{(\alpha_k - \alpha_m) t} \left( X_k(t) e^{-\alpha_k t} - \sum_{i=1}^{k-1} \frac{P_i^k}{a_1^i} W_i e^{(\alpha_i - \alpha_k) t} - W_k + \sum_{i=1}^{k-1} \frac{P_i^k}{a_1^i} W_i e^{(\alpha_i - \alpha_k) t} + W_k \right).$$

*On rearranging the terms, then* $X_m(t)e^{-\alpha_m t} + \sum_{i=1}^{m-1} a_i^m X_i(t)e^{-\alpha_m t} - W_m$ *equals*

$$= X_m(t)e^{-\alpha_m t} + a_1^m W_1 e^{(\alpha_1-\alpha_m)t} + \sum_{k=2}^{m-1} a_k^m \left( \sum_{i=1}^{k-1} \frac{P_i^k}{a_1^i} W_i e^{(\alpha_i-\alpha_m)t} + W_k e^{(\alpha_k-\alpha_m)t} \right) - W_m$$

$$+ a_1^m e^{(\alpha_1-\alpha_m)t} \left( X_1(t) e^{-\alpha_1 t} - W_1 \right) + \sum_{k=2}^{m-1} a_k^m e^{(\alpha_k-\alpha_m)t} \left( X_k(t) e^{-\alpha_k t} - \sum_{i=1}^{k-1} \frac{P_i^k}{a_1^i} W_i e^{(\alpha_i-\alpha_k)t} - W_k \right)$$

*Now using the following inequality: if $D = A+B+C$ then $|D| \geq |A|-|B|-|C|$ where $A, B, C, D$ are real numbers. In what follows, we can write*

$$\left| X_m(t)e^{-\alpha_m t} + \sum_{i=1}^{m-1} a_i^m X_i(t) e^{-\alpha_m t} - W_m \right|$$

$$\geq \left| X_m(t) e^{-\alpha_m t} + a_1^m W_1 e^{(\alpha_1-\alpha_m)t} + \sum_{k=2}^{m-1} a_k^m \left( \sum_{i=1}^{k-1} \frac{P_i^k}{a_1^i} W_i e^{(\alpha_i-\alpha_m)t} + W_k e^{(\alpha_k-\alpha_m)t} \right) - W_m \right|$$

$$- \left| a_1^m e^{(\alpha_1-\alpha_m)t} \left( X_1(t) e^{-\alpha_1 t} - W_1 \right) \right| - \sum_{k=2}^{m-1} \left| \left( a_k^m e^{(\alpha_k-\alpha_m)t} \right) \left( X_k(t) e^{-\alpha_k t} - \sum_{i=1}^{k-1} \frac{P_i^k}{a_1^i} W_i e^{(\alpha_i-\alpha_k)t} - W_k \right) \right|$$

$$= \left| X_m(t) e^{-\alpha_m t} + W_1 e^{(\alpha_1-\alpha_m)t} \left( a_1^m + \sum_{k=2}^{m-1} a_k^m P_1^k \right) + \sum_{k=3}^{m-1} \sum_{i=2}^{k-1} a_k^m \frac{P_i^k}{a_1^i} W_i e^{(\alpha_i-\alpha_m)t} + \sum_{k=2}^{m-1} W_k a_k^m e^{(\alpha_k-\alpha_m)t} - W_m \right|$$

$$- \left| a_1^m e^{(\alpha_1-\alpha_m)t} \left( X_1(t) e^{-\alpha_1 t} - W_1 \right) \right| - \sum_{k=2}^{m-1} \left| \left( a_k^m e^{(\alpha_k-\alpha_m)t} \right) \left( X_k(t) e^{-\alpha_k t} - \sum_{i=1}^{k-1} \frac{P_i^k}{a_1^i} W_i e^{(\alpha_i-\alpha_k)t} - W_k \right) \right|.$$

*Observe that* $a_1^m + \sum_{k=2}^{m-1} a_k^m P_1^k = -P_1^m$. *Changing the order of summations, we obtain*

$$\left| X_m(t) e^{-\alpha_m t} + \sum_{i=1}^{m-1} a_i^m X_i(t) e^{-\alpha_m t} - W_m \right|$$

$$\geq \left| X_m(t) - P_1^m W_1 e^{\alpha_1 t} + \sum_{i=2}^{m-2} \sum_{k=i+1}^{m-1} a_k^m \frac{P_i^k}{a_1^i} W_i e^{(\alpha_i-\alpha_m)t} + \sum_{k=2}^{m-1} W_k a_k^m e^{(\alpha_k-\alpha_m)t} - W_m \right|$$

$$- \left| a_1^m e^{(\alpha_1-\alpha_m)t} \left( X_1(t) e^{-\alpha_1 t} - W_1 \right) \right| - \sum_{k=2}^{m-1} \left| \left( a_k^m e^{(\alpha_k-\alpha_m)t} \right) \left( X_k(t) e^{-\alpha_k t} - \sum_{i=1}^{k-1} \frac{P_i^k}{a_1^i} W_i e^{(\alpha_i-\alpha_k)t} - W_k \right) \right|. \quad (45)$$

*Now observe the following*

$$\sum_{m=j+1}^{n-1} P_j^m a_m^n = \sum_{m=j+1}^{n-1} \sum_{k=1}^{m-j-1} (-1)^k \sum_{j_1=j+1}^{m-k} \cdots \sum_{j_k=j_{k-1}+1}^{m-1} a_1^j a_j^{j_1} a_{j_1}^{j_2} \cdots a_{j_k}^m a_m^n - \sum_{m=j+1}^{n-1} a_1^j a_j^m a_m^n.$$

*By changing a series of order of summations as before (e.g., see 40), we have*

$$\sum_{m=j+1}^{n-1} P_j^m a_m^n = \sum_{k=1}^{n-j-1} (-1)^k \sum_{j_1=j+1}^{n-k} \cdots \sum_{j_k=j_{k-1}+1}^{m-1} \sum_{m=j_k+1}^{n-1} a_1^j a_j^{j_1} a_{j_1}^{j_2} \cdots a_{j_k}^m a_m^n - \sum_{m=j+1}^{n-1} a_1^j a_j^m a_m^n. \quad (46)$$

*Thus,* $\sum_{m=j+1}^{n-1} P_j^m a_m^n + a_j^m a_m^n = -P_j^n$. *Substituting* $\sum_{m=j+1}^{n-1} P_j^m a_m^n = -P_j^n - a_j^m a_m^n$ *in the inequality as given in (45):*

$$\left| X_m(t) e^{-\alpha_m t} + \sum_{i=1}^{m-1} a_i^m X_i(t) e^{-\alpha_m t} - W_m \right|$$

$$\geq \left| X_m(t) e^{-\alpha_m t} - \sum_{i=1}^{m-1} \frac{P_i^m + a_j^m a_m^n}{a_1^i} W_i e^{(\alpha_i-\alpha_m)t} + \sum_{k=2}^{m-1} W_k a_k^m e^{(\alpha_k-\alpha_m)t} - W_m \right|$$

$$- \left| a_1^m e^{(\alpha_1-\alpha_m)t} \left( X_1(t) e^{-\alpha_1 t} - W_1 \right) \right| - \sum_{k=2}^{m-1} \left| \left( a_k^m e^{(\alpha_k-\alpha_m)t} \right) \left( X_k(t) e^{-\alpha_k t} - \sum_{i=1}^{k-1} \frac{P_i^k}{a_1^i} W_i e^{(\alpha_i-\alpha_k)t} - W_k \right) \right| \quad (47)$$

*Note that*

$$X_m(t)e^{-\alpha_m t} - \sum_{i=1}^{m-1} \frac{P_i^m + a_j^m a_m^n}{a_1^i} W_i e^{(\alpha_i - \alpha_m)t} + \sum_{k=2}^{m-1} W_k a_k^m e^{(\alpha_k - \alpha_m)t} - W_m$$

$$= \underbrace{X_m(t)e^{-\alpha_m t} - \sum_{i=1}^{m-1} \frac{P_i^m}{a_1^i} W_i e^{(\alpha_i - \alpha_m)t} - W_m}_{} + \underbrace{\sum_{k=2}^{m-1} W_k a_k^m e^{(\alpha_k - \alpha_m)t} - \sum_{i=1}^{m-1} \frac{a_j^m a_m^n}{a_1^i} W_i e^{(\alpha_i - \alpha_m)t}}_{}; \text{ as before}$$

$$\left| X_m(t)e^{-\alpha_m t} - \sum_{i=1}^{m-1} \frac{P_i^m + a_j^m a_m^n}{a_1^i} W_i e^{(\alpha_i - \alpha_m)t} + \sum_{k=2}^{m-1} W_k a_k^m e^{(\alpha_k - \alpha_m)t} - W_m \right|$$

$$\geq \left| X_m(t)e^{-\alpha_m t} - \sum_{i=1}^{m-1} \frac{P_i^m}{a_1^i} W_i e^{(\alpha_i - \alpha_m)t} - W_m \right| - \left| \sum_{k=2}^{m-1} W_k a_k^m e^{(\alpha_k - \alpha_m)t} \right| - \left| \sum_{i=1}^{m-1} \frac{a_j^m a_m^n}{a_1^i} W_i e^{(\alpha_i - \alpha_m)t} \right|$$

*In what follows, the inequality in the equation 47 becomes*

$$\left| X_m(t)e^{-\alpha_m t} + \sum_{i=1}^{m-1} a_i^m X_i(t) e^{-\alpha_m t} - W_m \right| \geq \left| X_m(t)e^{-\alpha_m t} - \sum_{i=1}^{m-1} \frac{P_i^m}{a_1^i} W_i e^{(\alpha_i - \alpha_m)t} - W_m \right|$$

$$- \left| \sum_{k=2}^{m-1} W_k a_k^m e^{(\alpha_k - \alpha_m)t} \right| - \left| \sum_{i=1}^{m-1} \frac{a_j^m a_m^n}{a_1^i} W_i e^{(\alpha_i - \alpha_m)t} \right|$$

$$- \left| a_1^m e^{(\alpha_1 - \alpha_m)t} \left( X_1(t)e^{-\alpha_1 t} - W_1 \right) \right| - \sum_{k=2}^{m-1} \left| \left( a_k^m e^{(\alpha_k - \alpha_m)t} \right) \left( X_k(t)e^{-\alpha_k t} - \sum_{i=1}^{k-1} \frac{P_i^k}{a_1^i} W_i e^{(\alpha_i - \alpha_k)t} - W_k \right) \right|$$

*Observe that the above inequality holds for all the sample paths, i.e., for any $\omega \in \Omega$, we have*

$$\left| X_m(t,\omega)e^{-\alpha_m t} + \sum_{i=1}^{m-1} a_i^m X_i(t,\omega) e^{-\alpha_m t} - W_m(\omega) \right| \geq \left| X_m(t,\omega)e^{-\alpha_m t} - \sum_{i=1}^{m-1} \frac{P_i^m}{a_1^i} W_i(\omega) e^{(\alpha_i - \alpha_m)t} - W_m(\omega) \right|$$

$$- \left| \sum_{k=2}^{m-1} W_k(\omega) a_k^m e^{(\alpha_k - \alpha_m)t} \right| - \left| \sum_{i=1}^{m-1} \frac{a_j^m a_m^n}{a_1^i} W_i(\omega) e^{(\alpha_i - \alpha_m)t} \right| - \left| a_1^m e^{(\alpha_1 - \alpha_m)t} \left( X_1(t,\omega)e^{-\alpha_1 t} - W_1(\omega) \right) \right|$$

$$- \sum_{k=2}^{m-1} \left| \left( a_k^m e^{(\alpha_k - \alpha_m)t} \right) \left( X_k(t,\omega)e^{-\alpha_k t} - \sum_{i=1}^{k-1} \frac{P_i^k}{a_1^i} W_i(\omega) e^{(\alpha_i - \alpha_k)t} - W_k(\omega) \right) \right|; \text{ thus}$$

$$\left| X_m(t)e^{-\alpha_m t} - \sum_{i=1}^{m-1} \frac{P_i^m}{a_1^i} W_i e^{(\alpha_i - \alpha_m)t} - W_m \right| \leq \left| X_m(t)e^{-\alpha_m t} + \sum_{i=1}^{m-1} a_i^m X_i(t) e^{-\alpha_m t} - W_m \right|$$

$$+ \left| \sum_{k=2}^{m-1} W_k(\omega) a_k^m e^{(\alpha_k - \alpha_m)t} \right| + \left| \sum_{i=1}^{m-1} \frac{a_j^m a_m^n}{a_1^i} W_i(\omega) e^{(\alpha_i - \alpha_m)t} \right| + \left| a_1^m e^{(\alpha_1 - \alpha_m)t} \left( X_1(t,\omega)e^{-\alpha_1 t} - W_1(\omega) \right) \right|$$

$$+ \sum_{k=2}^{m-1} \left| \left( a_k^m e^{(\alpha_k - \alpha_m)t} \right) \left( X_k(t,\omega)e^{-\alpha_k t} + \sum_{i=1}^{k-1} \frac{P_i^k}{a_1^i} W_i(\omega) e^{(\alpha_i - \alpha_k)t} - W_k(\omega) \right) \right|. \tag{48}$$

Appealing to the continuous mapping Theorem[10] in almost sure convergence (see Theorem 17.5 [14]), and using $\lim_{t \to \infty} e^{(\alpha_i - \alpha_k)t} = 0$ ($\because \alpha_k > \alpha_i \; \forall k > i$).

$$\left| a_1^m e^{(\alpha_1 - \alpha_m)t} \left( X_1(t)e^{-\alpha_1 t} - W_1 \right) \right| \xrightarrow{a.s.} 0 \text{ and}$$

$$\left| \left( a_k^m e^{(\alpha_k - \alpha_m)t} \right) \left( X_k(t)e^{-\alpha_k t} - \sum_{i=1}^{k-1} \frac{P_i^k}{a_1^i} W_i e^{(\alpha_i - \alpha_k)t} - W_k \right) \right| \xrightarrow{a.s.} 0 \; \forall \; k < m \text{ and so on.}$$

---

[10] If $X_n \xrightarrow{a.s.} X$, then $f(X_n) \xrightarrow{a.s.} f(X)$ where $f(.)$ is a continuous function.

Basically, all the terms on the right hand side of the equation 48 are converging to zero almost surely. Thus, we get the desired result

$$\left| X_m(t) e^{-\alpha_m t} - \sum_{i=1}^{m-1} \frac{P_i^m}{a_1^i} W_i e^{(\alpha_i - \alpha_m)t} - W_m \right| \xrightarrow{a.s.} 0.$$

And, for sufficiently large $t$, we have $\boxed{X_m(t) \approx \sum_{i=1}^{m-1} \frac{P_i^m}{a_1^i} W_i e^{\alpha_i t} + W_m e^{\alpha_m t}}.$

∎

**Lemma 1** *With $B_m$ as the generator matrix of the process for $2 \le m \le n$ and $\left(B_m^k\right)_{ji}$ be the $(j,i)$-th entry of the matrix $(B_m)^n$. We have $\left(e^{B_m \delta}\right)_{kk} = e^{\alpha_k \delta}$ $\forall$ $k$ and for $1 \le j < i \le m$*

$$\left(B_m^n\right)_{ji} = \alpha_{ji} \frac{\alpha_j^n - \alpha_i^n}{\alpha_j - \alpha_i} + \sum_{k=1}^{i-1-j} \sum_{j_1=j+1}^{i-k} \cdots \sum_{j_k=j_{k-1}+1}^{i-1} \sum_{p \in S_j^i(\mathbb{J}_k)} \frac{\alpha_{jj_1} \alpha_{j_1 j_2} \cdots \alpha_{j_k i}}{\prod_{l \in S_j^i(\mathbb{J}_k): l \ne p}(\alpha_p - \alpha_l)} \alpha_p^n; \text{ and}$$

$$\left(e^{B_m \delta}\right)_{ji} = \alpha_{ji} \frac{e^{\alpha_j \delta} - e^{\alpha_i \delta}}{\alpha_j - \alpha_i} + \sum_{k=1}^{i-1-j} \sum_{j_1=j+1}^{i-k} \cdots \sum_{j_k=j_{k-1}+1}^{i-1} \sum_{p \in S_j^i(\mathbb{J}_k)} \frac{\alpha_{jj_1} \alpha_{j_1 j_2} \cdots \alpha_{j_k i}}{\prod_{l \in S_j^i(\mathbb{J}_k): l \ne p}(\alpha_p - \alpha_l)} e^{\alpha_p \delta}. \quad (49)$$

**Proof 6** *We prove it through the method of induction. Firstly, note that one can decompose the generator matrix as below.*

$$B_m = \begin{bmatrix} \alpha_1 & \alpha_{12} & \cdots & \alpha_{1m-1} & \alpha_{1m} \\ 0 & \alpha_2 & \cdots & \alpha_{2m-1} & \alpha_{2m} \\ \vdots & \vdots & \ddots & \vdots & \vdots \\ 0 & 0 & \cdots & \alpha_{m-1} & \alpha_{m-1m} \\ \hline 0 & 0 & \cdots & 0 & \alpha_m \end{bmatrix} = \left[ \begin{array}{c|c} B_{m-1} & C \\ \hline \mathbf{0} & \alpha_m \end{array} \right]; \; C = [\alpha_{1m}, \alpha_{2m}, \cdots, \alpha_{m-1m}]. \quad (50)$$

*We require to compute $B_m^n$. With the help above decomposition (50), we have the following*

$$B_m^2 = \left[ \begin{array}{c|c} B_{m-1} & C \\ \hline \mathbf{0} & \alpha_m \end{array} \right] \left[ \begin{array}{c|c} B_{m-1} & C \\ \hline \mathbf{0} & \alpha_m \end{array} \right] = \left[ \begin{array}{c|c} B_{m-1}^2 & B_{m-1} C + \alpha_m C \\ \hline \mathbf{0} & \alpha_m^2 \end{array} \right]; \; \text{and in general}$$

$$B_m^n = \left[ \begin{array}{c|c} B_{m-1}^n & \left( B_{m-1}^{n-1} + \alpha_m B_{m-1}^{n-2} + \cdots + \alpha_m^{k-1} B_{m-1}^{n-k} + \cdots + \alpha_m^{n-2} B_{m-1} + \alpha_m^{n-1} I \right) C \\ \hline \mathbf{0} & \alpha_m^n \end{array} \right]; \quad (51)$$

*where $I$ is the identity matrix of the appropriate order.*
**Induction:** *For $m = 2$, we have $B_1 = \alpha_1, C = \alpha_{12}, \alpha_m = \alpha_2$ and*

$$B_2^n = \begin{bmatrix} \alpha_1^n & \alpha_{12} \frac{\alpha_1^n - \alpha_2^n}{\alpha_1 - \alpha_2} \\ 0 & \alpha_2^n \end{bmatrix}; \; \text{and } e^{B_2 \delta} = \begin{bmatrix} e^{\alpha_1 \delta} & \alpha_{12} \frac{e^{\alpha_1 \delta} - e^{\alpha_2 \delta}}{\alpha_1 - \alpha_2} \\ 0 & e^{\alpha_2 \delta} \end{bmatrix}. \quad (52)$$

*We see that each entry of $e^{B_2 \delta}$ is according to (49). Thus, the claim hold for $m = 2$.*

*Assume that result holds for $m - 1$, by which we have $\left(e^{B \delta}\right)_{kk} = e^{\alpha_k \delta}$ $\forall$ $k \le m-1$ and for $1 \le j < i \le m-1$*

$$\left(B_{m-1}^n\right)_{ji} = \alpha_{ji} \frac{\alpha_j^n - \alpha_i^n}{\alpha_j - \alpha_i} + \sum_{k=1}^{i-1-j} \sum_{j_1=j+1}^{i-k} \cdots \sum_{j_k=j_{k-1}+1}^{i-1} \sum_{p \in S_j^i(\mathbb{J}_k)} \frac{\alpha_{jj_1} \alpha_{j_1 j_2} \cdots \alpha_{j_k i} \alpha_j^n}{\prod_{l \in S_j^i(\mathbb{J}_k): l \ne p}(\alpha_p - \alpha_l)}; \text{ and} \quad (53)$$

$$\left(e^{B_{m-1} \delta}\right)_{ji} = \alpha_{ji} \frac{e^{\alpha_j \delta} - e^{\alpha_i \delta}}{\alpha_j - \alpha_i} + \sum_{k=1}^{i-1-j} \sum_{j_1=j+1}^{i-k} \cdots \sum_{j_k=j_{k-1}+1}^{i-1} \sum_{p \in S_j^i(\mathbb{J}_k)} \frac{\alpha_{jj_1} \alpha_{j_1 j_2} \cdots \alpha_{j_k i} e^{\alpha_j \delta}}{\prod_{l \in S_j^i(\mathbb{J}_k): l \ne p}(\alpha_p - \alpha_l)}. \quad (54)$$

*We have to prove it for $m$. From (51) it is clear that to compute $B_m^n$, we need to compute its last column only (as the rest of the entries are given by $B_{m-1}^n$). Further, observe that in the most general case, we require* to compute the first element of the last column of $e^{B_m\delta}$ only while the rest of the entries are filled by suitably changing the role of the indices. Utilizing both the facts, it suffices to compute $(B_m^n)_{1m}$ only in order to get the full matrix $B_m^n$ and hence $e^{B_m\delta}$. We thus calculate the first entry of the last column as follows

$$\left(B_{m-1}^{n-1} + \alpha_m B_{m-1}^{n-2} + \cdots + \alpha_m^{k-1} B_{m-1}^{n-k} + \cdots + \alpha_m^{n-2} B_{m-1} + \alpha_m^{n-1} I\right)C = \sum_{k_1=1}^{n} \alpha_m^{k_1-1} B_{m-1}^{n-k_1} C.$$

Note that it is obtained by multiplying the first row, say $[b_1, b_2, \cdots, b_{m-1}]$, of the matrix $\sum_{k_1=1}^{n} \alpha_m^{k_1-1} B_{m-1}^{n-k_1}$ with the vector $C$. Using equation (53) for $j = 1$, we have

$$b_1 = \sum_{k_1=1}^{n} \alpha_m^{k_1-1} \alpha_1^{n-k_1};$$

$$b_i = \sum_{k_1=1}^{n} \alpha_m^{k_1-1}\left(\alpha_{1i}\frac{\alpha_1^{n-k_1} - \alpha_i^{n-k_1}}{\alpha_1 - \alpha_i} + \sum_{k=1}^{i-2}\sum_{j_1=2}^{i-k}\cdots\sum_{j_k=j_{k-1}+1}^{i-1}\sum_{p\in S_1^i(\mathbb{J}_k)} \frac{\alpha_{1j_1}\alpha_{j_1 j_2}\cdots\alpha_{j_k i}\alpha_p^{n-k_1}}{\prod_{l\in S_1^i(\mathbb{J}_k): l\neq p}(\alpha_p - \alpha_l)}\right) \forall\, i \geq 2$$

$$= \frac{\alpha_{1i}\alpha_1^n}{(\alpha_1 - \alpha_i)(\alpha_1 - \alpha_m)} + \frac{\alpha_{1i}\alpha_i^n}{(\alpha_i - \alpha_1)(\alpha_i - \alpha_m)} + \frac{\alpha_{1i}\alpha_m^n}{(\alpha_m - \alpha_1)(\alpha_m - \alpha_i)}$$

$$+ \sum_{k=1}^{i-2}\sum_{j_1=2}^{i-k}\cdots\sum_{j_k=j_{k-1}+1}^{i-1}\sum_{p\in S_1^i(\mathbb{J}_k)} \frac{\alpha_{1j_1}\alpha_{j_1 j_2}\cdots\alpha_{j_k i}}{\prod_{l\in S_1^i(\mathbb{J}_k): l\neq p}(\alpha_p - \alpha_l)} \frac{\alpha_p^n - \alpha_m^n}{\alpha_p - \alpha_m}.$$

And

$$(B_m^n)_{1m} = \sum_{i=1}^{m-1} b_i \alpha_{im} = \alpha_{1m}\frac{\alpha_1^n - \alpha_m^n}{\alpha_1 - \alpha_m} + \sum_{i=2}^{m-1}\left(\frac{\alpha_{1i}\alpha_{im}\alpha_1^n}{(\alpha_1 - \alpha_i)(\alpha_1 - \alpha_m)} + \frac{\alpha_{1i}\alpha_{im}\alpha_i^n}{(\alpha_i - \alpha_1)(\alpha_i - \alpha_m)} + \frac{\alpha_{1i}\alpha_{im}\alpha_m^n}{(\alpha_m - \alpha_1)(\alpha_m - \alpha_i)}\right)$$

$$+ \sum_{i=3}^{m-1}\sum_{k=1}^{i-2}\sum_{j_1=2}^{i-k}\cdots\sum_{j_k=j_{k-1}+1}^{i-1}\sum_{p\in S_1^i(\mathbb{J}_k)} \frac{\alpha_{1j_1}\alpha_{j_1 j_2}\cdots\alpha_{j_k i}\alpha_{im}}{\prod_{l\in S_1^i(\mathbb{J}_k): l\neq p}(\alpha_p - \alpha_l)} \frac{\alpha_p^n - \alpha_m^n}{\alpha_p - \alpha_m}$$

$$= \alpha_{1m}\frac{\alpha_1^n - \alpha_m^n}{\alpha_1 - \alpha_m} + \sum_{k=1}^{m-2}\sum_{j_1=2}^{m-k}\sum_{j_2=j_1+1}^{m-k+1}\cdots\sum_{j_k=j_{k-1}+1}^{m-1}\sum_{p\in S_1^m(\mathbb{J}_k)} \frac{\alpha_{1j_1}\alpha_{j_1 j_2}\cdots\alpha_{j_k m}\alpha_p^n}{\prod_{l\in S_1^m(\mathbb{J}_k): l\neq p}(\alpha_p - \alpha_l)}.$$

Observe that

$$\sum_{j\in S_1^i(\mathbb{J}_k)} \frac{1}{\prod_{l\in S_1^i(\mathbb{J}_k): i\neq j}(\alpha_j - \alpha_l)} = 0 \,\forall\, i \text{ and } \mathbb{J}_k \neq \emptyset. \tag{55}$$

It is easy to verify the following using the expression for $(B_m^n)_{1m}$ and the equation (55)

$$\left(e^{B_m\delta}\right)_{1m} = \alpha_{1m}\frac{e^{\alpha_1} - e^{\alpha_m}}{\alpha_1 - \alpha_m} + \sum_{k=1}^{m-2}\sum_{j_1=2}^{m-k}\sum_{j_2=j_1+1}^{m-k+1}\cdots\sum_{j_k=j_{k-1}+1}^{m-1}\sum_{p\in S_1^m(\mathbb{J}_k)} \frac{\alpha_{1j_1}\alpha_{j_1 j_2}\cdots\alpha_{j_k m}e^{\alpha_p}}{\prod_{l\in S_1^m(\mathbb{J}_k): l\neq p}(\alpha_p - \alpha_l)}.$$

We see that $\left(e^{B_m\delta}\right)_{1m}$ matches the entry as per the hypothesis. And the remaining entries of the last column of $e^{B_m\delta}$ are computed using the following observation.

$$B_m = \begin{bmatrix} \alpha_1 & \alpha_{12} & \alpha_{13} & \cdots & \alpha_{1m-1} & \alpha_{1m} \\ \hline 0 & \alpha_2 & \alpha_{23} & \cdots & \alpha_{2m-1} & \alpha_{2m} \\ 0 & 0 & \alpha_3 & \cdots & \alpha_{3m-1} & \alpha_{3m} \\ \vdots & \vdots & & \ddots & \vdots & \vdots \\ 0 & 0 & 0 & \cdots & \alpha_{m-1} & \alpha_{m-1m} \\ 0 & 0 & 0 & \cdots & 0 & \alpha_m \end{bmatrix}.$$

The lower $m-1 \times m-1$ block of the matrix $B_m$ is similar to $B_{m-1}$ when all the row and column indices are shifted above by one each. Further, as the matrices, i.e. $B_{m-1}, B_m$, are upper triangular, the lower $m-1 \times m-1$ block of $e^{B_m}$ preserves the structure of the matrix $e^{B_{m-1}}$ with the said indices shifting. In all, we have for $1 \leq j < i \leq m$

$$\left(e^{B_m\delta}\right)_{ji} = \alpha_{ji}\frac{e^{\alpha_j\delta} - e^{\alpha_i\delta}}{\alpha_j - \alpha_i} + \sum_{k=1}^{i-1-j}\sum_{j_1=j+1}^{i-k}\cdots\sum_{j_k=j_{k-1}+1}^{i-1}\sum_{p\in S_j^i(\mathbb{J}_k)} \frac{\alpha_{jj_1}\alpha_{j_1 j_2}\cdots\alpha_{j_k i}}{\prod_{l\in S_j^i(\mathbb{J}_k): l\neq p}(\alpha_p - \alpha_l)} e^{\alpha_p\delta}. \tag{56}$$

*Hence proved.*

■

**Proof 7** *of Theorem 2: This proof is constructed using the standard arguments used for the extinction probability in branching processes. Appealing to Theorem 1, we have*

$$\left| X_m(t) e^{-\alpha_m t} + \sum_{i=1}^{m-1} a_i^m X_i(t) e^{-\alpha_m t} - W_m \right| \xrightarrow{a.s.} 0; \text{ and} \tag{57}$$

$$\left| X_m(t) e^{-\alpha_m t} - \sum_{i=1}^{m-1} \frac{P_i^m}{a_1^i} W_i e^{(\alpha_i - \alpha_m)t} - W_m \right| \xrightarrow{a.s.} 0 \tag{58}$$

*where $W_1, W_2, \cdots, W_m$ are non-negative random variables. As we know that the above result holds when the process is started by a type-1 particle (most general framework). We have similar results when the process is initiated by particle(s) of some other type. For the sake of this proof, we require the equivalent of Theorem 1 when the starting particle is of type-i. Towards this, let $X_{ij}(t)$ represent the size of the type-j population at time t, when the process starts with one type-i particle. Then the equivalent of (58) is:*

$$\left| X_{ij}(t) e^{-\alpha_j t} - \sum_{k=i}^{j-1} a_k^j X_{kj} e^{-\alpha_j t} - W_{ij} \right| \xrightarrow{a.s.} 0$$

*and*

$$\left| X_{ij}(t) e^{-\alpha_j t} - \sum_{k=i}^{j-1} \frac{P_k^j}{a_i^k} W_{kj} e^{(\alpha_k - \alpha_j)t} - W_{ij} \right| \xrightarrow{a.s.} 0.$$

*where $W_{ij}, W_{i+1 j}, \cdots, W_{jj}$ are non-negative and integrable random variables. Observe that*

$$X_{ij}(t) = 0 \text{ for some } t > 0 \iff W_{ij} = W_{i+1 j} = \cdots = W_{jj} = 0. \tag{59}$$

*This is because $\{e^{\alpha_j}\}_{1 \leq j \leq m}$ are linearly independent, and $P_j^m, a_j^m$ are non zero. Further, as $\alpha_1 < \alpha_2 \cdots < \alpha_m$, it is easy to see that*

$$\lim_{t \to \infty} X_{ij}(t) e^{-\alpha_j t} = W_{ij}. \tag{60}$$

*As we have from equation (3)*

$$q_{ij} = P[X_i(t) = 0 \ \forall \ i \leq j \text{ for some } t > 0 \mid \mathbf{X}(0) = \mathbf{e}_i]. \tag{61}$$

*Recall that $X_{mm}(t)$ evolves according to a standard single-type CTBP. As we have $E[\zeta_{ij} \log \zeta_{ij}] < \infty \ \forall \ i, j$, using the well-known result for characterization of the random variable $W_m$ in irreducible BP ( [2]), we have $P(W_{mm} > 0) = q_{mm}$. And we precisely show that an equivalent of that holds in our framework too.*

*Say $\zeta_{ij}$ be the number of type-j offspring produced by one type-i particle ($j \geq i$) at the first transition epoch $\tau$, then one can write:*

$$X_{ij}(\tau + t) = \sum_{k=i}^{j} \sum_{l=1}^{\zeta_{ik}} X_{kj}^l(t);$$

*where $X_{kj}^l(t)$ be the size of type-j produced by the l-th particle of type-k*

$$X_{ij}(\tau+t) e^{-\alpha_j(\tau+t)} = \sum_{k=i}^{j} \sum_{l=1}^{\zeta_{ik}} X_{kj}^l(t) e^{-\alpha_j(\tau+t)}; \text{ now taking the limit on both the sides}$$

$$\lim_{t \to \infty} X_{ij}(\tau+t) e^{-\alpha_j(\tau+t)} = \lim_{t \to \infty} \sum_{k=i}^{j} \sum_{l=1}^{\zeta_{ik}} X_{kj}^l(t) e^{-\alpha_j(\tau+t)}.$$

*And using equation (60),* $W_{ij} = \sum_{k=i}^{j} \sum_{l=1}^{\zeta_{ik}} W_{kj}^l e^{-\alpha_j \tau};$ \hfill (62)

*where $\{W_{kj}^l\}_l$, for any given k and j, are IID random variables.*

Observe that $W_{ij} = 0$ if and only if $W_{ij}^l = 0, W_{i+1j}^l = 0, \cdots, W_{jj}^l = 0$ for all $l$. Define the following

$$\rho_{ij} := P(W_{ij}(\omega) = 0) = P\left(\cap_{k=i}^{j} \cap_{l=1}^{\zeta_i} W_{kj}^l(\omega) = 0\right). \tag{63}$$

One can write $\rho_{ij}$ by conditioning on the offspring of various types which are independent of each other

$$\begin{aligned}\rho_{ij} &= \sum_{k_{i+1},\cdots k_j} P\big(W_{ij}(\omega) = 0 \big| \zeta_{ii+1} = k_{i+1}, \zeta_{ii+2} = k_{i+2}\cdots, \zeta_{ij} = k_j\big) P\big(\zeta_{ii+1} = k_{i+1},\cdots, \zeta_{ij} = k_j\big) \\ &= \sum_{k_{i+1},\cdots k_j} \rho_{ii+1}^{k_{i+1}} \rho_{ii+2}^{k_{i+2}} \cdots \rho_{ij}^{k_j} P(\zeta_{ii+1} = k_{i+1}, \zeta_{ii+2} = k_{i+2}\cdots, \zeta_{ij} = k_j)\end{aligned}$$

which is the same equation satisfied by the extinction probabilities (4). Hence proved. ∎

**Proof 8** *of Theorem 3: Using the linear property of expectation*

$$\begin{aligned}E\left[Y(t+\delta)^T \xi_2^R + X(t+\delta)^T (\alpha_2 I - A_{11})^{-1} A_{12}\xi_2^R \big| X(t), Y(t)\right] &= E\left[Y(t+\delta)^T \xi_2^R \big| X(t), Y(t)\right] \\ &+ E\left[X(t+\delta)^T (\alpha_2 I - A_{11})^{-1} A_{12}\xi_2^R \big| X(t), Y(t)\right]\end{aligned} \tag{64}$$

*and*

$$\begin{aligned}E\left[Y(t+\delta)^T \xi_2^R \big| X(t), Y(t)\right] &= Y(t)^T e^{A_{22}\delta} \xi_2^R + E\left[Y(t+\delta)^T \xi_2^R \big| X(t)\right] \\ &= Y(t)^T \xi_2^R e^{\alpha_2 \delta} + X(t)^T (e^{A\delta})_{12} \xi_2^R \\ E\left[X(t+\delta)^T (\alpha_2 I - A_{11})^{-1} A_{12}\xi_2^R \big| X(t), Y(t)\right] &= E\left[X(t+\delta)\big|X(t), Y(t)\right]^T (\alpha_2 I - A_{11})^{-1} A_{12}\xi_2^R \\ &= X(t)^T e^{A_{11}\delta} (\alpha_2 I - A_{11})^{-1} A_{12}\xi_2^R.\end{aligned}$$

*As we know that $e^{\mathbb{A}} = I + \mathbb{A} + \frac{\mathbb{A}^2}{2!} + \frac{\mathbb{A}^3}{3!} + \cdots$. We now compute $\mathbb{A}^n$ as follows*

$$\mathbb{A}^2 = \begin{bmatrix} A_{11} & A_{12} \\ \mathbf{0} & A_{22} \end{bmatrix} \begin{bmatrix} A_{11} & A_{12} \\ \mathbf{0} & A_{22} \end{bmatrix} = \begin{bmatrix} A_{11}^2 & A_{11}A_{12} + A_{12}A_{22} \\ \mathbf{0} & A_{22}^2 \end{bmatrix}$$

*and $\mathbb{A}_{12}^3 = A_{11}^2 A_{12} + A_{11}A_{12}A_{22} + A_{12}A_{22}^2$ so on. Using $A_{22}^n \xi_2^R = \alpha_2^n \xi_2^R$, we get the following*

$$\begin{aligned}(e^A)_{12}\xi_2^R &= \frac{A_{12}\xi_2^R}{1!} + \frac{A_{11}A_{12}\xi_2^R + A_{12}\alpha_2\xi_2^R}{2!} + \frac{A_{11}^2 A_{12}\xi_2^R + A_{11}A_{12}\alpha_2\xi_2^R + A_{12}\alpha_2^2\xi_2^R}{3!} + \cdots \\ &= \left[I + \frac{A_{11} + \alpha_2 I}{2!} + \frac{A_{11}^2 + A_{11}\alpha_2 + \alpha_2^2 I}{3!} + \cdots\right] A_{12}\xi_2^R \\ &= \frac{(I-B_{11})^{-1}}{\alpha_2}\left[\alpha_2(I-B_{11}) + \frac{\alpha_2^2}{2!}\left(I-B_{11}^2\right) + \frac{\alpha_2^3}{3!}\left(I-B_{11}^3\right) + \frac{\alpha_2^4}{4!}\left(I-B_{11}^4\right) + \cdots\right] A_{12}\xi_2^R \quad (B_{11} = A_{11}/\alpha_2) \\ &= (\alpha_2 I - A_{11})^{-1}\left(e^{\alpha_2}I - e^{A_{11}}\right) A_{12}\xi_2^R \\ (e^{\mathbb{A}\delta})_{12}\xi_2 &= (\alpha_2 I - A_{11})^{-1}\left(e^{\alpha_2 \delta}I - e^{A_{11}\delta}\right) A_{12}\xi_2^R.\end{aligned}$$

*On substituting all these in (64) we get*

$$\begin{aligned}e^{-\alpha_2(t+\delta)} E\left[Y(t+\delta)^T \xi_2^R + X(t+\delta)^T (\alpha_2 I - A_{11})^{-1} A_{12}\xi_2^R \big| X(t), Y(t)\right] \\ = e^{-\alpha_2(t+\delta)} X(t)^T (\alpha_2 I - A_{11})^{-1}\left(e^{\alpha_2 \delta}I - e^{A_{11}\delta}\right) A_{12}\xi_2^R \\ + e^{-\alpha_2(t+\delta)}\left(Y(t)^T \xi_2^R e^{\alpha_2 \delta} + X(t)^T e^{A_{11}\delta}(\alpha_2 I - A_{11})^{-1} A_{12}\xi_2^R\right) \\ = X(t)^T (\alpha_2 I - A_{11})^{-1} e^{-\alpha_2 t} A_{12}\xi_2^R + e^{-\alpha_2 t} Y(t)^T \xi_2^R.\end{aligned} \tag{65}$$

*This completes the martingale part.*
**Convergence:** *When $\alpha_2 > \alpha_1$, then the matrix $(\alpha_2 I - A_{11})$ becomes an M-matrix whose inverse is a non-negative matrix. This is because $(\alpha_2 I - A_{11})$ is a Z-matrix and the real part of all its eigen values are strictly* positive (definition in page 3 and Lemma 4.1 in [13]). Thus matrix $(\alpha_2 I - A_{11})^{-1}$ has nonnegative entries, $A_{12}$ has nonnegative entries and by Perron Frobenius theory $\xi_2$ is a positive vector and thus (65) is a nonnegative vector. Thus, we have a non-negative martingale and by forward martingale convergence theorem it converges to $\tilde{W}_{12}$ which is non-negative and integrable. Thus,

$$\left| e^{-\alpha_2 t} \boldsymbol{Y}(t)\xi_2^R + e^{-\alpha_2 t} \boldsymbol{X}(t)(\alpha_2 I - A_{11})^{-1} A_{12}\xi_2^R - \tilde{W}_{12} \right| \xrightarrow{a.s.} 0. \tag{66}$$

As before, we have

$$\left| e^{-\alpha_2 t} \boldsymbol{Y}(t)\xi_2^R + e^{-\alpha_2 t} W_1 e^{\alpha_1 t} \xi_1^R (\alpha_2 I - A_{11})^{-1} A_{12}\xi_2^R - \tilde{W}_{12} \right|$$

$$= \left| e^{-\alpha_2 t} \boldsymbol{Y}(t)\xi_2^R + e^{-\alpha_2 t} \left( \boldsymbol{X}(t) - \boldsymbol{X}(t) + W_1 e^{\alpha_1 t} \xi_1^R \right) (\alpha_2 I - A_{11})^{-1} A_{12}\xi_2^R - \tilde{W}_{12} \right|$$

$$\leq \left| e^{-\alpha_2 t} \boldsymbol{Y}(t)\xi_2^R + e^{-\alpha_2 t} \boldsymbol{X}(t)(\alpha_2 I - A_{11})^{-1} A_{12}\xi_2^R - \tilde{W}_{12} \right|$$

$$+ \left| e^{-\alpha_2 t} \left( -\boldsymbol{X}(t) + W_1 e^{\alpha_1 t} \xi_1^R \right) (\alpha_2 I - A_{11})^{-1} A_{12}\xi_2^R \right|.$$

As $t \to \infty$, we obtain (as before) $\left| e^{-\alpha_2 t} \boldsymbol{Y}(t)\xi_2^R + e^{-\alpha_2 t} W_1 e^{\alpha_1 t} \xi_1^R (\alpha_2 I - A_{11})^{-1} A_{12}\xi_2^R - \tilde{W}_{12} \right| \xrightarrow{a.s.} 0$. Equivalently, it is easy to see that

$$\left| \boldsymbol{Y}(t)\xi_2^R + W_1 e^{\alpha_1 t} \xi_1^R (\alpha_2 I - A_{11})^{-1} A_{12}\xi_2^R - \tilde{W}_{12} e^{\alpha_2 t} \right| \xrightarrow{a.s.} 0.$$

Further, for sufficiently large $t$, we have the following approximation

$$\boldsymbol{Y}(t)^T \xi_2^R \approx \tilde{W}_{12} e^{\alpha_2 t} - W_1 e^{\alpha_1 t} \xi_1^R (\alpha_2 I - A_{11})^{-1} A_{12}\xi_2^R.$$

■